\documentclass[12pt]{amsart}

\usepackage{amsfonts,amsmath,amssymb,amsthm}
\usepackage{latexsym}
\usepackage{mathrsfs}
\usepackage{tikz}
\usepackage{color}
\usepackage{enumitem}

\newcommand{\gsl}{\mathfrak{sl}}

\newcommand{\Usl}{U_q(\gsl_2)}
\newcommand{\Utrp}{U_q(\gsl_2)^{\otimes 3}}
\newcommand{\C}{\Gamma}
\newcommand{\M}[1]{M_{#1}}

\newcommand{\cR}{\mathcal{R}}

\newcommand{\cS}{\mathcal{S}}
\newcommand{\cK}{\mathcal{K}}

\newcommand{\cJ}{\mathcal{J}}
\newcommand{\cM}{\mathcal{M}}
\newcommand{\cC}{\mathcal{C}}

\newcommand{\mC}{\mathfrak{C}}
\newcommand{\gC}{\mathfrak{C}}

\newcommand{\fa}{\mathfrak{a}}
\newcommand{\fb}{\mathfrak{b}}
\newcommand{\fc}{\mathfrak{c}}

\newcommand{\tA}{\tilde{A}}
\newcommand{\tB}{\tilde{B}}
\newcommand{\tD}{\tilde{D}}
\newcommand{\tal}{\tilde{\alpha}}
\newcommand{\tK}{\tilde{K}}
\newcommand{\tG}{\tilde{G}}

\newcommand{\tchi}{\tilde{\chi}}

\newcommand{\qi}{{q^{-1}}}
\newcommand{\qnum}[1]{\frac{q^{#1}-q^{-#1}}{q-q^{-1}}}

\newcommand{\kpt}[3]{#1\otimes #2 \otimes #3}
\newcommand{\kpd}[2]{#1\otimes #2}

\newcommand{\id}{\text{id}}
\newcommand{\idr}[1]{\mathbb{I}_{2#1+1}}

\newcommand{\End}{\text{End}}

\newtheorem{prop}{Proposition}[section]
\newtheorem{thm}{Theorem}[section]
\newtheorem{rem}{Remark}[section]

\newtheorem{definition}{Definition}[section]
\newtheorem{conj}{Conjecture}[section]

\numberwithin{equation}{section}

\begin{document}

\title[Askey--Wilson algebra and centralizers]{Temperley--Lieb, Birman--Murakami--Wenzl and Askey--Wilson algebras and other centralizers of $\Usl$}

\author[N.Cramp\'e]{Nicolas Cramp\'e$^{\dagger}$}
\address{$^\dagger$ Institut Denis-Poisson CNRS/UMR 7013 - Universit\'e de Tours - Universit\'e d'Orl\'eans, 
	Parc de Grandmont, 37200 Tours, France.}
\email{crampe1977@gmail.com}

\author[L.Vinet]{Luc Vinet$^{*}$}
\address{$^*$ Centre de recherches math\'ematiques, Universit\'e de Montr\'eal,
	P.O. Box 6128, Centre-ville Station,
	Montr\'eal (Qu\'ebec), H3C 3J7, Canada.}
\email{vinet@crm.umontreal.ca}

\author[M.Zaimi]{Meri Zaimi$^{**}$}
\address{$^{**}$ Centre de recherches math\'ematiques, Universit\'e de Montr\'eal,
	P.O. Box 6128, Centre-ville Station,
	Montr\'eal (Qu\'ebec), H3C 3J7, Canada.}
\email{meri.zaimi@umontreal.ca}

\begin{abstract}
	The centralizer of the image of the diagonal embedding of $U_q(\gsl_2)$ in the tensor product of three irreducible representations is examined in a Schur--Weyl duality spirit. The aim is to offer a description in terms of generators and relations. A conjecture in this respect is offered with the centralizers presented as quotients of the Askey--Wilson algebra. Support for the conjecture is provided by an examination of the representations of the quotients. The conjecture is also shown to be true in a number of cases thereby exhibiting in particular the Temperley--Lieb, Birman--Murakami--Wenzl and one-boundary Temperley--Lieb algebras as quotients of the Askey--Wilson algebra.	
\end{abstract}

\maketitle

\vspace{3mm}

\section{Introduction}

The objective of this paper is to establish precisely the connections between the Askey--Wilson algebra and the centralizers of the quantum algebra $U_q(\gsl_2)$ such as the Temperley--Lieb and Birman--Murakami--Wenzl algebras. 

In previous works, the connections between the Racah, Temperley--Lieb and Brauer algebras and other centralizers of $\gsl_2$ were studied in the spirit of the Schur--Weyl duality \cite{CPV}. In a similar fashion, the Bannai-Ito algebra was connected to the centralizers of the superalgebra $\mathfrak{osp}(1|2)$, and in particular to the Brauer algebra \cite{CFV}. The present paper generalizes the results of \cite{CPV} by examining their $q$-deformation.  

The Askey--Wislon algebra was first introduced in \cite{Zh} and is defined by three generators satisfying some $q$-commutation relations. This algebra encodes the properties of the Askey--Wilson polynomials \cite{Koek} and is related to the Racah problem for $U_q(\gsl_2)$ \cite{GZ}. Due to this connection, a centrally extended Askey--Wilson algebra can be mapped to the centralizer of the diagonal embedding of $U_q(\gsl_2)$ into $\Utrp$ \cite{Ter,H1}. In the $q$-deformation of the universal enveloping algebra of $\gsl_2$ to the quantum algebra $U_q(\gsl_2)$, the  Askey--Wilson algebra plays a role analogous to that of the Racah algebra.

From the Schur--Weyl duality, the centralizer of the diagonal embedding of $U_q(\gsl_2)$ in the tensor product of its fundamental representation is connected to the Hecke algebra. In the case of the threefold tensor product, the centralizer is known \cite{Jimbo} to be isomorphic to the Temperley--Lieb algebra \cite{TL}, which is a quotient of the Hecke algebra. In fact, the algebra $U_q(\gsl_2)$ has infinitely many finite irreducible representations, labeled by a half-integer or integer spin $j$. In the case of the tensor product of three spin-$1$ representations, it is also known \cite{LZ} that the centralizer of $U_q(\gsl_2)$ is isomorphic to the Birman--Murakami--Wenzl algebra \cite{IMO}, which is a $q$-deformation of the Brauer algebra. However, in the general case of three irreducible representations of spins $j_1,j_2$ and $j_3$, an algebraic description of the centralizer is not known.

The present paper provides an attempt to describe the centralizer of the image of the diagonal embedding of $U_q(\gsl_2)$ in the tensor product of any three irreducible representations in terms of generators and relations, by using the connections with the Askey--Wilson algebra. It is first shown that there is a surjective map (between generators) from the Askey--Wilson algebra to the centralizer. This statement corresponds in invariant theory to the first fundamental theorem \cite{KP}. A conjecture is then proposed in order to obtain an isomorphism between a quotient (given in terms of relations) of the Askey--Wilson algebra and the centralizer of $U_q(\gsl_2)$ -- this relates to the second fundamental theorem in invariant theory \cite{KP}. The conjecture is proved for three spin-$\frac{1}{2}$ representations, in which case the Temperley--Lieb algebra is obtained explicitly as a quotient of the Askey--Wilson algebra. Similarly, for three spin-$1$ representations, it is shown that the conjecture holds and that the Birman--Murakami--Wenzl algebra is isomorphic to a quotient of the Askey--Wilson algebra. The conjecture is also verified for three spin-$\frac{3}{2}$ representations, and for one spin-$j$ and two spin-$\frac{1}{2}$ representations, for $j$ any spin greater than $\frac{1}{2}$. In the latter case, it is shown that the centralizer is isomorphic to the one-boundary Tempereley--Lieb algebra \cite{MS,MW,NRG}.  

The plan of this paper is as follows. Section \ref{sec:UslAW} gives the precise connection between the centralizer of $\Usl$ and the Askey--Wilson algebra. Subsection \ref{ssec:Usl} presents the quantum algebra $\Usl$ and its properties. The centralizer of $\Usl$ in $\Utrp$ and the intermediate Casimirs are defined in Subsection \ref{ssec:CUsl}. A homomorphism between the centrally extended Askey--Wilson algebra $AW(3)$ and this centralizer is given in Subsection \ref{ssec:CAW}.
Section \ref{sec:IrrepC} is concerned with the representations of elements in $\Utrp$. The finite irreducible representations of $\Usl$ and their tensor product decomposition rules are recalled in Subsection \ref{ssec:irreps}. Subsection \ref{ssec:centr} introduces the object of main interest, that is the centralizer of the image of the diagonal embedding of $\Usl$ in the tensor product of three irreducible representations. 
Section \ref{sec:conj} aims to describe this centralizer in terms of generators and relations. Subsection \ref{ssec:homo} maps $AW(3)$ to this centralizer, and this is shown to be a surjection in Subsection \ref{ssec:surj}. The kernel of this map is discussed in Subsection \ref{ssec:ker}, and a conjecture proposing that a quotient of $AW(3)$ is isomorphic to the centralizer is formulated. Subsection \ref{ssec:perm} contains the proof that the conjecture does not depend on the ordering of the three spins $j_1,j_2,j_3$. In order to support the conjecture, Section \ref{sec:repAW} studies the finite irreducible representations of the quotient of $AW(3)$. The remaining sections contain the proofs of the conjecture for some particular cases. Section \ref{sec:AW1TL} focuses on the case $j_1=j_2=j_3=\frac{1}{2}$. It is shown in Subsection \ref{ssec:AW1} that the conjecture holds in this case, and the precise connection with the Temperley--Lieb algebra is given in Subsection \ref{ssec:TL}. Section \ref{sec:AW2BMW} considers the case $j_1=j_2=j_3=1$. The proof of the conjecture is given in Subsection \ref{ssec:AW2}. An isomorphism between the quotient of $AW(3)$ and the Birman--Murakami--Wenzl algebra is obtained in Subsection \ref{ssec:BMW}. The conjecture for the case $j_1=j_2=j_3=\frac{3}{2}$ is proved in Section \ref{sec:AW3}. Finally, Section \ref{sec:AWj1bTL} studies the case $j_1=j$ for $j=1,\frac{3}{2},...$ and $j_2=j_3=\frac{1}{2}$. The conjecture is verified in Subsection \ref{ssec:AWj}, and the connection with the one-boundary Temperley--Lieb algebra is described in Subsection \ref{ssec:1bTL}.        

\section{Centralizer of $U_q(\gsl_2)$ and Askey--Wilson algebra}\label{sec:UslAW}

In this section, we recall well-known properties of the quantum algebra $U_q(\gsl_2)$ to fix the notations.
Then, the definition of the centralizer of the diagonal embedding of $U_q(\gsl_2)$ in  $U_q(\gsl_2)^{\otimes 3}$ is recalled
and its homomorphism with the centrally extended Askey--Wilson algebra $AW(3)$ is presented.

\subsection{$U_q(\gsl_2)$ algebra}\label{ssec:Usl}

The associative algebra $U_q(\gsl_2)$ is generated by $E$, $F$ and $q^{H}$ with the defining relations
\begin{equation}
q^{H}E=q  Eq^H \ , \quad q^{H}F=q^{-1} Fq^H  \ \text{and} \quad [E,F]=[2H]_q \ , \label{eq:Uqsl2}
\end{equation}
where $[X]_q=\qnum{X}$. Throughout this paper, $q$ is a complex number not root of unity.
There is a central element in $U_q(\gsl_2)$, called quadratic Casimir element, given by
\begin{equation}
\C=(q-\qi)^2 FE+qq^{2H}+\qi q^{-2H} \ . \label{eq:C}
\end{equation}

There exists also an algebra homomorphism $\Delta:U_q(\gsl_2) \rightarrow U_q(\gsl_2)\otimes U_q(\gsl_2) $, called comultiplication, defined on the generators by 
\begin{equation}
\Delta(E)=E \otimes q^{-H} +q^H \otimes E \ , \quad  \Delta(F)=F \otimes q^{-H} +q^{H} \otimes F \quad \text{and} \quad \Delta(q^H)=q^H\otimes q^H\ . \label{eq:com1}
\end{equation}
This comultiplication is coassociative 
\begin{equation}
(\Delta \otimes \id)\Delta= (\id \otimes \Delta)\Delta =: \Delta^{(2)} \ .\label{eq:com2}
\end{equation}

We define the opposite comultiplication $\Delta^{op}=\sigma \circ \Delta$, where $\sigma(x \otimes y)=y \otimes x$, for $x,y\in \Usl$.
It is a homomorphism from $U_q(\gsl_2)$ to $U_q(\gsl_2)\otimes U_q(\gsl_2)$ different from $\Delta$. 
Both are related by the universal $R$-matrix $\cR\in U_q(\gsl_2) \otimes U_q(\gsl_2)$ satisfying
\begin{equation}
\Delta(x) \cR  = \cR \Delta^{op}(x) \quad \text{for } x\in U_q(\gsl_2)\label{eq:RD} \ .
\end{equation}
The universal $R$-matrix also satisfies the Yang-Baxter equation
\begin{equation}
 \cR_{12}\cR_{13}\cR_{23}=\cR_{23}\cR_{13}\cR_{12}\ . \label{eq:YBE}
\end{equation}
We have used the usual notations: if $\cR=\cR^\alpha\otimes \cR_\alpha$, then $\cR_{12}=\cR^\alpha\otimes \cR_\alpha \otimes 1$, $\cR_{23}=1\otimes \cR^\alpha\otimes \cR_\alpha$
and $\cR_{13}=\cR^\alpha\otimes1\otimes  \cR_\alpha$ (the sum w.r.t. $\alpha$ is understood).

\subsection{Centralizer of $\Usl$ in $\Utrp$}\label{ssec:CUsl}

The centralizer $\gC_3$ of the diagonal embedding of $U_q(\gsl_2)$ in $U_q(\gsl_2)^{\otimes 3}$ is
\begin{equation}
\gC_3=\{ X \in U_q(\gsl_2)^{\otimes 3}\ \big| \ [\Delta^{(2)}(x), X]=0\ , \ \  \forall x\in U_q(\gsl_2) \} . \label{eq:centr1}
\end{equation}
This centralizer is a subalgebra of $\Utrp$ and we want to describe this subalgebra with some generators and defining relations.
Let us first give some elements of $\gC_3$ by using the Casimir element $\C$ which is central in $\Usl$.
We define the following Casimir elements of $\Utrp$
\begin{equation}
\C_1=\C\otimes 1 \otimes 1 \ , \qquad \C_2=1\otimes \C \otimes 1 \ , \qquad \C_3=1\otimes 1 \otimes \C \ .
\label{eq:Ci}
\end{equation}
These elements are central in $\Utrp$ and thus belong to $\gC_3$. We also define the total Casimir
\begin{equation}
\C_{123}=\Delta^{(2)}(\C) \ . \label{eq:C123}
\end{equation}
This element belongs to $\gC_3$ because
$[\Delta^{(2)}(\C),\Delta^{(2)}(x)]=\Delta^{(2)}([\C,x])=0$ for all $x\in U_q(\gsl_2)$. Let us notice that $\C_{123}$ is central in $\mC_3$ since it is also an element of the diagonal embedding of $\Usl$.

We then define the intermediate Casimirs associated to the recoupling of the two first or the two last factors of $\Utrp$
\begin{equation}
\C_{12}=\Delta(\C) \otimes 1  \quad \text{and} \qquad \C_{23}=1 \otimes \Delta(\C) \ . \label{eq:Cint}
\end{equation}
One uses the properties of the comultiplication to show that $\C_{12}$ and $\C_{23}$ are in $\gC_3$; indeed, for all $x\in U_q(\gsl_2)$,
\begin{align}
[\C_{12},\Delta^{(2)}(x)]=[\Delta(\C) \otimes 1,(\Delta \otimes \id)\Delta(x)]=
(\Delta \otimes \id)[\C\otimes 1,\Delta(x)]=0 \ \label{eq:C12ctz},\\
[\C_{23},\Delta^{(2)}(x)]=
[1 \otimes \Delta(\C),(\id \otimes \Delta)\Delta(x)]=
(\id \otimes \Delta)[1\otimes \C,\Delta(x)]=0 \ . \label{eq:C23ctz}
\end{align}
In the limit $q\to1$, it can be shown that the element
\begin{equation}
\C_{13}=\sum_{\alpha}\C_{\alpha}\otimes 1 \otimes \C^{\alpha} \ , \label{eq:C13}
\end{equation}
where $\Delta(\C)=\sum_{\alpha} \C_{\alpha}\otimes \C^{\alpha}$, 
belongs to the centralizer of $U(\gsl_2)$ in $U(\gsl_2)^{\otimes 3}$.
However, this is not the case for the quantum algebra $\Usl$. This difficulty that arises in the $q$-deformation of the algebra $U(\gsl_2)$
was addressed in \cite{CGVZ} where a definition of the third intermediate Casimir element of $\Usl$ is provided with the help of the universal $R$-matrix. It is shown in \cite{CGVZ} that the following elements 
 \begin{eqnarray}
\C_{13}^{(0)}&=\cR_{12}  \C_{13} \cR_{12}^{-1}=\cR_{32}^{-1}  \C_{13} \cR_{32}\ ,\label{eq:C130}\\
\C_{13}^{(1)}&=\cR_{23}  \C_{13} \cR_{23}^{-1}=\cR_{21}^{-1}  \C_{13} \cR_{21}\ ,\label{eq:C131}
\end{eqnarray}
are in the centralizer $\gC_3$.

\subsection{Connection with the Askey--Wilson algebra}\label{ssec:CAW} 

The intermediate Casimir elements $\C_{12}$, $\C_{23}$, $\C_{13}^{(0)}$ and $\C_{13}^{(1)}$ do not commute pairwise but satisfy certain relations
which are identified as those of the Askey--Wilson algebra $AW(3)$.
\begin{definition}\label{def:aw3}
	The centrally extended Askey--Wilson algebra $AW(3)$ is generated by $A$, $B$, $D$ and central elements $\alpha_1$, $\alpha_2$, $\alpha_3$ and $K$ subject to the following defining relations
	\begin{align}
	A +\frac{[B,D]_q}{q^2-q^{-2}} &=\frac{\alpha_1\alpha_2+\alpha_3K}{q+q^{-1}} \ , \label{eq:AWcBD}\\
	B +\frac{[D,A]_q}{q^2-q^{-2}} &=\frac{\alpha_2\alpha_3+\alpha_1K}{q+q^{-1}} \ , \label{eq:AWcDA} \\
	D +\frac{[A,B]_q}{q^2-q^{-2}} &=\frac{\alpha_1\alpha_3+\alpha_2 K}{q+q^{-1}} \ , \label{eq:AWcAB}
	\end{align}
	where $[X,Y]_q=qXY-\qi YX$. We also define the element $D'\in AW(3)$ by the following relation
	\begin{align}
	D' +\frac{[B,A]_q}{q^2-q^{-2}} &=\frac{\alpha_1\alpha_3+\alpha_2 K}{q+q^{-1}}\ . \label{eq:AWcBA}
	\end{align}
\end{definition}
The algebra $AW(3)$ has a Casimir element given by
\begin{equation}
	\Omega = qA (\alpha_1\alpha_2+\alpha_3 K) +q^{-1} B (\alpha_2\alpha_3+\alpha_1 K) + q D(\alpha_1\alpha_3+\alpha_2 K) -q^2 A^2 - q^{-2} B^2 -q^2 D^2 - q ABD \ . \label{eq:Om}
\end{equation}

The connection between the centralizer $\gC_3$ defined by \eqref{eq:centr1} and the Askey--Wilson algebra is given in the following proposition.
\begin{prop}\label{pr:varphi} The map 
	$\varphi:AW(3)\rightarrow \gC_3$ defined by
	\begin{align}
	\varphi(\alpha_i)=\C_i \ , \  \varphi(A)=\C_{12} \ , \  \varphi(B)=\C_{23} \ , \  \varphi(K)=\C_{123} \ , \label{eq:homoAWC}
	\end{align}
	is an algebra homomorphism. We deduce that
	\begin{align}
 \varphi(D)=\C_{13}^{(0)} \ , \  \varphi(D')=\C_{13}^{(1)} \ . \  \label{eq:DR}
	\end{align}
\end{prop}
\noindent  The homomorphism has been proved in \cite{GZ}; a direct computation shows that the intermediate Casimir elements satisfy all the relations of $AW(3)$.
Relations \eqref{eq:DR} have been proved more recently in \cite{CGVZ}
and a simpler proof of the homomorphism using the universal $R$-matrix has also been  given. 
Let us remark that a similar proof has also been simplified in the case of the Bannai--Ito algebra 
and the centralizer for the super Lie algebra $\mathfrak{osp}(1|2)$ \cite{CVZ}.

Using \eqref{eq:AWcAB} to replace $D$ in \eqref{eq:AWcBD} and \eqref{eq:AWcDA}, one shows that the following relations provide an equivalent presentation of $AW(3)$ which will be useful for later computations
\begin{align}
\frac{[B,[A,B]_q]_q}{(q-\qi)^2}=&(q+\qi)^2A+(\alpha_1 \alpha_3 +\alpha_2 K)B-(q+\qi)(\alpha_1\alpha_2+\alpha_3 K) \ , \label{eq:AWd1}\\
\frac{[[A,B]_q,A]_q}{(q-\qi)^2}=&(q+\qi)^2B+(\alpha_1 \alpha_3 +\alpha_2 K)A-(q+\qi)(\alpha_2\alpha_3+\alpha_1 K) \ . \label{eq:AWd2}
\end{align}
Furthermore, noticing that $[X,[Y,X]_q]_q=[[X,Y]_q,X]_q$ and using the element $D'$ defined in \eqref{eq:AWcBA}, one finds that \eqref{eq:AWd1} and \eqref{eq:AWd2} imply
\begin{align}
A +\frac{[D',B]_q}{q^2-q^{-2}} &=\frac{\alpha_1\alpha_2+\alpha_3K}{q+q^{-1}} \ , \label{eq:AWcD'B}\\
B +\frac{[A,D']_q}{q^2-q^{-2}} &=\frac{\alpha_2\alpha_3+\alpha_1K}{q+q^{-1}} \ . \label{eq:AWcAD'}
\end{align}
Relations \eqref{eq:AWcBA}, \eqref{eq:AWcD'B} and \eqref{eq:AWcAD'} provide another $\mathbb{Z}_3$ symmetric presentation of $AW(3)$.

\begin{rem}\label{rem:qlim}
	Upon performing the affine transformation $X=(q-\qi)^2 \tilde{X}+q+\qi$ on the elements $X=A,B,D,D',\alpha_i,K$ of $AW(3)$, one sees that relations \eqref{eq:AWd1}--\eqref{eq:AWd2} can be written as
	\begin{align}
	[\tilde{B},[\tA,\tilde{B}]_q]_q=&(q+\qi)\left(-\tilde{B}^2-\{\tA,\tilde{B}\}+(\tK+\tal_1+\tal_2+\tal_3)\tilde{B}+(\tal_1-\tK)(\tal_3-\tal_2)\right) \label{eq:AWtd1} \\
	+&(q-\qi)^2(\tal_1\tal_3+\tal_2\tK)\tilde{B} \ ,  \nonumber\\
	[[\tA,\tilde{B}]_q,\tA]_q=&(q+\qi)\left(-\tA^2-\{\tA,\tilde{B}\}+(\tK+\tal_1+\tal_2+\tal_3)\tA+(\tal_3-\tK)(\tal_1-\tal_2)\right) \label{eq:AWtd2}\\
	+&(q-\qi)^2(\tal_1\tal_3+\tal_2\tK)\tA \ , \nonumber
	\end{align}
	where $\{X,Y\}=XY+YX$. By taking the limit $q\to1$ of \eqref{eq:AWtd1} and \eqref{eq:AWtd2}, one recovers the defining relations of the Racah algebra used in \cite{CPV}.
	Relations \eqref{eq:AWcAB} and \eqref{eq:AWcBA} are transformed into
	\begin{align}
	\frac{[\tA,\tB]_q}{q-\qi}&=\tal_1\tal_3+\tal_2\tK+\frac{q+\qi}{(q-\qi)^2}(\tal_1+\tal_2+\tal_3+\tK-\tA-\tB-\tD) \ , \label{eq:AWtcAB}\\
	\frac{[\tB,\tA]_q}{q-\qi}&=\tal_1\tal_3+\tal_2\tK+\frac{q+\qi}{(q-\qi)^2}(\tal_1+\tal_2+\tal_3+\tK-\tA-\tB-\tD')\ . \label{eq:AWtcBA}
	\end{align}
	In the limit $q\to1$, the elements $\tD$ and $\tD'$ are equal, and the images by $\varphi$ of \eqref{eq:AWtcAB} and \eqref{eq:AWtcBA} both reduce 
	to the well-known linear relation $\tilde{\C}_1+\tilde{\C}_2+\tilde{\C}_3+\tilde{\C}_{123}-\tilde{\C}_{12}-\tilde{\C}_{23}-\tilde{\C}_{13}=0$ that holds in $U(\gsl_2)^{\otimes 3}$.
	 
\end{rem}

\section{Decomposition of tensor product of representations and centralizer}\label{sec:IrrepC}

In the previous section, we introduced the centralizer $\gC_3$ of the diagonal embedding of $\Usl$ in $\Utrp$ and 
showed the connection of this subalgebra of $\Utrp$ with the Askey--Wilson algebra $AW(3)$.
We now focus on the corresponding objects when each $\Usl$ in $\Utrp$ is taken in a finite irreducible representation.

\subsection{Finite irreducible representations of $\Usl$}\label{ssec:irreps}

The quantum algebra $\Usl$ has finite irreducible representations of dimension $2j+1$ that we will denote by $\M{j}$, with $j=0,\frac{1}{2},1,\frac{3}{2},...$ 
The representation map will be denoted by $\pi_{j}:U_q(\gsl_2)\rightarrow \text{End}(\M{j})$. We will use the name spin-$j$ representation to refer to $\M{j}$. The representation of the Casimir element \eqref{eq:C} in the space $\M{j}$ is
\begin{equation}
\pi_{j}(\C)= \chi_j \idr{j}\quad \text{where } \chi_j= q^{2j+1}+q^{-2j-1}, 
\label{eq:casimirrep}
\end{equation}
and $\idr{j}$ is the $2j+1$ by $2j+1$ identity matrix.
We define the following sets, for three half-integers or integers $j_1$, $j_2$ and $j_3$
\begin{eqnarray}
&&\cJ (j_1, j_2 )=  \{|j_1-j_2|,|j_1-j_2|+1,...,j_1+j_2\} \ , \label{eq:J12}\\
&&\cJ(j_1,j_2,j_3)=\bigcup_{j \in\cJ(j_1,j_2)}\cJ(j,j_3) \ . \label{eq:J123}
\end{eqnarray}
Notice that there are no repeated numbers in $\cJ(j_1,j_2,j_3)$, and this set is invariant under any permutation of $j_1$, $j_2$ and $j_3$.

For $q$ not a root of unity,
the tensor product of two irreducible representations of $\Usl$ decomposes into the following 
direct sum of irreducible representations
\begin{equation}
\M{j_1}\otimes\M{j_2}=\bigoplus_{j\in\cJ(j_1,j_2)}\M{j} \ . \label{eq:decomp2}
\end{equation}
Similarly, the threefold tensor product of irreducible representations of $\Usl$ decomposes into the following direct sum
\begin{equation}
\M{j_1}\otimes\M{j_2}\otimes\M{j_3}=\bigoplus_{j\in\cJ(j_1,j_2,j_3)}d_j\M{j} \ , \label{eq:decomp3}
\end{equation}
where $d_j\in\mathbb{Z}_{>0}$ is the degeneracy of $\M{j}$ and is referred to as the Littlewood--Richardson coefficient.

\subsection{Centralizer of $U_q(\gsl_2)$ in $\End(\kpt{\M{j_1}}{\M{j_2}}{\M{j_3}})$ \label{ssec:centr}}

From now on, we fix three half-integers or integers $j_1$, $j_2$ and $j_3$. The centralizer $\cC_{j_1,j_2,j_3}$ of the diagonal embedding of $U_q(\gsl_2)$ in $\End(\kpt{\M{j_1}}{\M{j_2}}{\M{j_3}})$ is
\begin{equation}
\cC_{j_1,j_2,j_3}=\{m \in \End(\kpt{\M{j_1}}{\M{j_2}}{\M{j_3}}) \  \big| \  [\pi_{j_1,j_2,j_3} (\Delta^{(2)}(x)),m]=0,\forall x\in U_q(\gsl_2)\} \ , \label{eq:Cent}
\end{equation}
where we have used the shortened notation $\pi_{j_1,j_2,j_3}=\kpt{\pi_{j_1}}{\pi_{j_2}}{\pi_{j_3}} $. 
This centralizer as a subalgebra of $\End(\kpt{\M{j_1}}{\M{j_2}}{\M{j_3}})$ is the object of interest of this paper. 
In the next section, we conjecture a presentation of this centralizer 
in terms of generators and relations by using the connections with the Askey--Wilson algebra $AW(3)$. 
Let us first recall some known properties of this centralizer.

The knowledge of the centralizer permits to write the decomposition rule \eqref{eq:decomp3} as follows
\begin{equation}
\M{j_1}\otimes\M{j_2}\otimes\M{j_3}=\bigoplus_{j\in\cJ(j_1,j_2,j_3)} \M{j} \otimes V_j \ , \label{eq:decomp4}
\end{equation}
where $V_j$ is a finite irreducible representation of dimension $d_j$ of $\cC_{j_1,j_2,j_3}$. 
The set $\{V_j \ |\ j\in\cJ(j_1,j_2,j_3)\}$ is the complete set of non-equivalent irreducible representations of $\cC_{j_1,j_2,j_3}$.
In particular, one deduces that the dimension of the centralizer is
\begin{equation}
\text{dim}(\cC_{j_1,j_2,j_3})=\sum_{j\in\cJ(j_1,j_2,j_3)}d_j^2 \ . \label{eq:dimCent}
\end{equation}
These representations $V_j$ are explicitly given in Subsection \ref{ssec:surj}.

We now define the images of the centralizing elements \eqref{eq:Ci}--\eqref{eq:Cint} and \eqref{eq:C130}--\eqref{eq:C131} of $\Utrp$ in the representation $\End(\kpt{\M{j_1}}{\M{j_2}}{\M{j_3}})$ as follows
\begin{equation}
	\begin{array}{rll}
	\pi_{j_1,j_2,j_3}:\hspace{1cm} \cC_3\hspace{1cm} &\rightarrow & \cC_{j_1,j_2,j_3} \\
	\Gamma_i, \Gamma_{ij}, \Gamma_{123} & \mapsto & C_i, C_{ij}, C_{123} \ .
	\end{array}
\end{equation}
Therefore, $\cC_{j_1,j_2,j_3}$ contains the elements $C_i$, $C_{ij}$ and $C_{123}$.
According to \eqref{eq:casimirrep}, the elements $C_i$ are simply constant 
matrices of value $\chi_{j_i}$ for $i=1,2,3$. The intermediate Casimirs 
$C_{12}$, $C_{23}$, $C_{13}^{(0)}$ and $C_{13}^{(1)}$ and the total Casimir $C_{123}$ of $\cC_{j_1,j_2,j_3}$ 
can be diagonalized if $q$ is not a root of unity.

Since $C_{12}$ is the Casimir associated to the recoupling of the two first factors of the 
threefold tensor product of $\Usl$, one finds (using the decomposition rule \eqref{eq:decomp2}) 
that its eigenvalues are $\chi_j$ for $j\in\cJ(j_1,j_2)$. Similarly, the eigenvalues of $C_{23}$ (resp. $C_{123}$)  
are $\chi_j$ for $j\in\cJ(j_2,j_3)$ (resp. $\cJ(j_1,j_2,j_3)$). 
The same argument cannot be applied directly to the intermediate Casimirs $C_{13}^{(0)}$ and $C_{13}^{(1)}$ since they are not 
trivial in the space 2. However, the element $C_{13}$ defined in \eqref{eq:C13} only couples the 
spaces 1 and 3 such that its eigenvalues are $\chi_j$ for $j\in\cJ(j_1,j_3)$. From the definitions \eqref{eq:C130} 
and \eqref{eq:C131}, we see that $C_{13}^{(0)}$ and $C_{13}^{(1)}$ are both conjugations of $C_{13}$ by an $R$-matrix. 
Hence, their eigenvalues are the same as those of $C_{13}$. 
The previous discussion implies that the minimal polynomials of the intermediate Casimirs and the total Casimir take the following form
\begin{align}
&\prod_{j\in\cJ(j_1,j_2)}(C_{12}-\chi_j)=0 \ , \quad \prod_{j\in\cJ(j_2,j_3}(C_{23}-\chi_j)=0 \ ,
\prod_{j\in\cJ(j_1,j_2,j_3)}(C_{123}-\chi_j)=0  \ , 
\label{eq:CH1}\\
&\prod_{j\in\cJ(j_1,j_3)}(C_{13}^{(0)}-\chi_j)=0  \ , \quad
\prod_{j\in\cJ(j_1,j_3)}(C_{13}^{(1)}-\chi_j)=0  \ . \label{eq:CH2}
\end{align}

Because $C_{123}$ is central in $\cC_{j_1,j_2,j_3}$, it can be diagonalized simultaneously with $C_{12}$, $C_{23}$, $C_{13}^{(0)}$ or $C_{13}^{(1)}$. Therefore, one gets the following minimal polynomials
\begin{align}
&\prod_{m\in\cM(j_1,j_2,j_3)}(C_{123}- C_{12}-m)=0 \ , \quad \prod_{m\in\cM(j_2,j_3,j_1)}(C_{123}- C_{23}-m)=0 \ ,\label{eq:CH3}\\
&\prod_{m\in\cM(j_1,j_3,j_2)}(C_{123}- C_{13}^{(0)}-m)=0 \ , \quad \prod_{m\in \cM(j_1,j_3,j_2)}(C_{123}- C_{13}^{(1)}-m)=0 \ , \label{eq:CH4}
\end{align}
where
\begin{equation}
\cM(j_a,j_b,j_c)=\bigcup_{j \in\cJ(j_a,j_b)}\{\chi_\ell-\chi_{j} \big| \ell \in \mathcal{J}(j ,j_c)\} \ . \label{eq:M123}
\end{equation}
In the previous set $\cM(j_a,j_b,j_c)$, there are no repeated numbers.

Before concluding this section, let us notice that if one performs the transformation given in Remark \ref{rem:qlim} on the Casimir element $\Gamma$ of $\Usl$, its value in the representation $\End(M_j)$ is $\tilde{\chi}_j=[j]_q[j+1]_q$. By construction, similar results hold for the eigenvalues of the transformed elements $\tilde{C}_i,\tilde{C}_{ij}$ and $\tilde{C}_{123}$. In the limit $q\to1$, the minimal polynomials of these transformed elements thus reduce to the ones discussed in \cite{CPV}.

\section{Algebraic description of the centralizer $\cC_{j_1,j_2,j_3}$ \label{sec:conj}}

Take $j_1$, $j_2$ and $j_3$ to be three fixed half-integers or integers. This section contains an attempt to give a definition of the centralizer $\cC_{j_1,j_2,j_3}$
in terms of generators and relations. We rely on the connection with the Askey--Wilson algebra $AW(3)$.

\subsection{Homomorphism with $AW(3)$}\label{ssec:homo}

The intermediate Casimir elements $C_i$, $C_{ij}$ and $C_{123}$ belonging to $\cC_{j_1,j_2,j_3}$ satisfy the 
defining relations of the Askey--Wilson algebra as stated precisely in the following proposition.
\begin{prop}\label{pr:phi} The map 
	$\phi:AW(3)\rightarrow \cC_{j_1,j_2,j_3}$ defined by
	\begin{align}
	\phi(\alpha_i)=C_i \ , \  \phi(A)=C_{12} \ , \  \phi(B)=C_{23} \ , \  \phi(K)=C_{123} \ , \label{eq:homoAWC3}
	\end{align}
	is an algebra homomorphism.
\end{prop}
\proof 
The result follows from the fact that $\phi$ is the composition of two homomorphisms $\phi=\pi_{j_1,j_2,j_3}\circ \varphi$, where $\varphi$ is defined in Proposition \ref{pr:varphi}.
\endproof
We recall that $C_i$ is $\chi_{j_i}=q^{2j_i+1}+q^{-2j_i-1}$ times the identity matrix. Therefore, it can be identified as the number $\chi_{j_i}$.   
Let us also emphasize that $\phi(D)=C_{13}^{(0)}$ and $\phi(D')=C_{13}^{(1)}$. Moreover, the image by $\phi$ of the Casimir element $\Omega$ of $AW(3)$ defined in \eqref{eq:Om} is equal to an expression involving only central elements \cite{GZ} :
\begin{equation}
\phi(\Omega)=C_1^2+C_2^2+C_3^2+C_{123}^2+C_1C_2C_3C_{123} - (q+\qi)^2. \label{eq:phiOm}
\end{equation}

\subsection{Surjectivity\label{ssec:surj}}

We now show that the intermediate Casimir elements $C_i$, $C_{12}$, $C_{23}$ and $C_{123}$ generate the whole centralizer $\cC_{j_1,j_2,j_3}$.
\begin{prop}\label{pr:surj} The map 
	$\phi:AW(3)\rightarrow \cC_{j_1,j_2,j_3}$ is surjective.
\end{prop}

\proof To reach that conclusion, we prove that the dimension of the image of $\phi$ is at least $\displaystyle \sum_{\ell\in \cJ(j_1,j_2,j_3)} d_\ell^2$, the dimension of $\cC_{j_1,j_2,j_3}$. 
Let $\ell\in \cJ(j_1,j_2,j_3)$ and 
\begin{equation}
 \cS^\ell(j_1,j_2,j_3) = \{\  j \in \cJ (j_1,  j_2 )\ |\ \ell \in \cJ(j,j_3) \ \}\ . \label{eq:defSk}
\end{equation}
From the definition \eqref{eq:J12}, we deduce that $\cS^\ell(j_1,j_2,j_3)=\{j_{\text{min}},j_{\text{min}}+1,\dots, j_{\text{max}} \}$ with
\begin{equation}
j_{\text{min}}= \text{max}( |j_1-j_2| , |j_3-\ell| )\quad \text{and} \qquad  j_{\text{max}}= \text{min}( j_1+j_2 , j_3+\ell ) \ . \label{eq:jminmax}
\end{equation}
The cardinality of this set is $d_\ell=j_{\text{max}}-j_{\text{min}}+1$. We denote by $M^+_\ell$ the vector space spanned by the highest weight vectors of the representations $M_\ell$ in the decomposition \eqref{eq:decomp3}.
The dimension of $M_\ell^+$ is $d_\ell$ and we can choose $d_\ell$ independent vectors $v_j\in M^+_\ell$ with $j\in \cS^\ell(j_1,j_2,j_3)$ such that 
\begin{equation}
\pi_{j_1,j_2,j_3} (\Delta^{(2)}(E))v_j=0\ , \quad \pi_{j_1,j_2,j_3} (\Delta^{(2)}(q^H))v_j=q^{\ell}v_j \ , \quad  C_{123} v_j=\chi_\ell v_j\ , \quad C_{12} v_j=\chi_j v_j \ ,
\end{equation}
and 
\begin{equation}
 C_{23} v_j = \sum_{k\in \cS^\ell(j_1,j_2,j_3)}  \alpha_{j,k} v_k \ ,
\end{equation}
where $\alpha_{j,k}$ are complex numbers.
The elements $C_{12}$ and $C_{23}$ are the images by $\phi$ of $A$ and $B$. Therefore, they satisfy the Askey--Wilson algebra. 
It is enough to determine the constants $\alpha_{j,k}$ as shown previously in \cite{Zh}. We reproduce this computation in the particular case needed here. We define the constants
\begin{equation}
	\fa=\chi_{j_1}\chi_{j_2}+\chi_{j_3}\chi_{\ell} \ , \quad \fb=\chi_{j_2}\chi_{j_3}+\chi_{j_1}\chi_{\ell} \quad \text{and} \quad \fc=\chi_{j_1}\chi_{j_3}+\chi_{j_2}\chi_{\ell} \ .
\end{equation}
We act with relation \eqref{eq:AWd2} on the vector $v_j$ (for $j \in \cS^\ell(j_1,j_2,j_3)$) and project the result on $v_k$ 
with $k\neq j$ and on $v_j$. We get 
\begin{eqnarray}
&&\Big( [j+k+2]_q[j+k]_q[k-j-1]_q[k-j+1]_q\Big) \alpha_{j,k}=0 \ , \label{eq:aljk}\\
&&  \alpha_{j,j}=\frac{\fc\chi_j-\fb\chi_0}{\chi_j^2-\chi_0^2} \qquad \text{for}\quad  j\neq 0 \ . \label{eq:aljj}
\end{eqnarray}
The projection on $v_j$ is trivial if $j=0$.
From relation \eqref{eq:aljk}, we deduce that $\alpha_{j,k}=0$ for $j \in \cS^\ell(j_1,j_2,j_3)$ and $k\neq  j+1,j-1,j$.

Then, we act with relation \eqref{eq:AWd1} on the vector $v_j$ and project the result on $v_{j-2}, v_{j-1}, \dots, v_{j+2}$. The projections are trivial except the one on $v_j$ which gives
the following relation
\begin{equation}
[2j+3]_q \alpha_{j,j+1}\alpha_{j+1,j}-[2j-1]_q \alpha_{j-1,j}\alpha_{j,j-1}= \frac{1}{\chi_0}(\fc - \chi_j \alpha_{j,j})\alpha_{j,j} +\chi_0\chi_j - \fa \ , \label{eq:arec}
\end{equation}
with the boundary conditions $\alpha_{j_{\text{min}},j_{\text{min}}-1}=0$ and $\alpha_{j_{\text{max}},j_{\text{max}}+1}=0$. 
By using \eqref{eq:aljj}, one can show that the recurrence relation \eqref{eq:arec} and the boundary condition $\alpha_{j_{\text{max}},j_{\text{max}}+1}=0$ imply
\begin{equation}
\alpha_{j-1,j}\alpha_{j,j-1}=\frac{\prod_{i=1}^{4}([j-r_i]_q[j+r_i]_q)}{[2j-1]_q[2j]_q^2[2j+1]_q}(q-\qi)^4 \qquad \text{for}\quad  j\neq 0 \ , \label{eq:aaj}
\end{equation}
where $r_1=j_1-j_2$, $r_2=j_3-\ell$, $r_3=j_1+j_2+1$ and $r_4=\ell+j_3+1$. We see from \eqref{eq:jminmax} that the second boundary condition $\alpha_{j_{\text{min}},j_{\text{min}}-1}=0$ is automatically satisfied
if $j_{\text{min}}>0$. In the case where $j_{\text{min}}=0$ (which only happens if $j_1=j_2$ and $\ell=j_3$), the limit $j\to0$ of \eqref{eq:aaj} vanishes. Moreover, we can deduce from \eqref{eq:arec} that $\alpha_{0,0}=\chi_{j_1}\chi_{j_3}/\chi_0$, which is the limit $j\to0$ of \eqref{eq:aljj}.

To conclude the proof, we notice that equation \eqref{eq:jminmax} implies that the R.H.S. of relation \eqref{eq:aaj} is never zero for $j_\text{min}<j\leq j_\text{max}$, 
and that the eigenvalues of $C_{12}$ are pairwise distinct.
Therefore, for a given $\ell \in \cJ(j_1,j_2,j_3)$, $C_{12}$ and $C_{23}$ generate a vector space of dimension $d_\ell^2$.
\endproof

\subsection{Kernel}\label{ssec:ker}

The map $\phi$ defined in the Proposition \ref{pr:phi} is not injective since there are non-trivial elements of $AW(3)$ 
that are mapped to zero, as seen from the results \eqref{eq:CH1}--\eqref{eq:CH4}. We want to provide 
a description of the kernel of the map $\phi$ in order to find a quotient of $AW(3)$ that is isomorphic to the centralizer $\cC_{j_1,j_2,j_3}$.
Let us first define a quotient of $AW(3)$.
\begin{definition}
	The algebra $\overline{AW}(j_1,j_2,j_3)$ is the quotient of the centrally extended Askey--Wilson algebra $AW(3)$ by the following relations
	\begin{align}
	&\alpha_i=\chi_{j_i} \ , \label{eq:AWq0} \\ 
	&\prod_{j\in\cJ(j_1,j_2)}(A-\chi_j)=0 \ , \prod_{j\in\cJ(j_2,j_3)}(B-\chi_j)=0 \ , \prod_{j\in\cJ(j_1,j_2,j_3)}(K-\chi_j)=0 \ , \label{eq:AWq1}\\
	&\prod_{j\in\cJ(j_1,j_3)}\left(D-\chi_j \right)=0 \ , \prod_{j\in\cJ(j_1,j_3)}\left(D'-\chi_j\right)=0 \ , \label{eq:AWq2}\\
	&\prod_{m\in\cM(j_1,j_2,j_3)}(K-A-m)=0 \ , \quad \prod_{m\in\cM(j_2,j_3,j_1)}(K-B-m)=0 \ , \label{eq:AWq3} \\
	&\prod_{m\in\cM(j_1,j_3,j_2)}\left(K-D-m\right)=0 \ ,\quad \prod_{m\in\cM(j_1,j_3,j_2)}\left(K-D'-m\right)=0 \ , \label{eq:AWq4} \\
	&\Omega=\chi_{j_1}^2+\chi_{j_2}^2+\chi_{j_3}^2+K^2+\chi_{j_1}\chi_{j_2}\chi_{j_3}K-(q+q^{-1})^2 \ , \label{eq:relOm}
	\end{align}
	where we recall that $D$ and $D'$ are defined through \eqref{eq:AWcAB}--\eqref{eq:AWcBA}, and $\Omega$ is defined in \eqref{eq:Om}.
\end{definition}
Let us emphasize that all the relations \eqref{eq:AWq0}--\eqref{eq:relOm} are in the kernel of the map $\phi$ in view of the results of Subsections \ref{ssec:centr} and \ref{ssec:homo}. We are now in position to state a conjecture that proposes an algebraic description of $\cC_{j_1,j_2,j_3}$.  
\begin{conj}\label{conj1}
The map $\overline{\phi}:\overline{AW}(j_1,j_2,j_3)\to \cC_{j_1,j_2,j_3}$ given by
	\begin{align}
 \overline{\phi}(A)=C_{12} \ , \  \overline{\phi}(B)=C_{23} \ , \  \overline{\phi}(K)=C_{123} \ , \label{eq:homoAWC2}
	\end{align}
is an algebra isomorphism.
\end{conj}
To support this conjecture, we remark that by taking the limit $q\to 1$ (as described in Remark \ref{rem:qlim}) of relations \eqref{eq:AWq0}--\eqref{eq:AWq4}, we recover the conjecture proposed in \cite{CPV} and proved in numerous cases. Let us notice that in this limit, the two relations in \eqref{eq:AWq2} reduce to only one relation, and similarly for the two relations in \eqref{eq:AWq4}. The relation \eqref{eq:relOm} involving the Casimir element of $AW(3)$ is new and will be useful when we discuss the representations of $\overline{AW}(j_1,j_2,j_3)$ in Section \ref{sec:repAW}. 

From the previous results, we know that $\overline{\phi}$ is a surjective homomorphism.
It remains to prove that it is injective, which can be done by demonstrating that
\begin{equation}
 \text{dim}(\overline{AW}(j_1,j_2,j_3)) \leq \sum_{j\in \cJ(j_1,j_2,j_3)} d_j^2=\text{dim}(\cC_{j_1,j_2,j_3}) \ . \label{eq:min1}
\end{equation}
To simplify the demonstration of \eqref{eq:min1}, we can decompose $\overline{AW}(j_1,j_2,j_3)$ into a direct sum of simpler algebras.
Indeed, let us introduce the following central idempotents, for $k \in \cJ(j_1,j_2,j_3)$,
\begin{equation}
 \cK_k=\prod_{\genfrac{}{}{0pt}{}{r \in \cJ(j_1,j_2,j_3)}{r\neq k}  } \frac{K -\chi_r}{\chi_k-\chi_r}\ ,
\end{equation}
which satisfy  $\cK_k\cK_\ell=\delta_{k,\ell}\cK_k$, $\displaystyle \sum_{k \in \cJ(j_1,j_2,j_3)} \cK_k=1$ and $K \cK_k=\cK_k K =\chi_k \cK_k$.
We deduce that 
\begin{equation}
  \overline{AW}(j_1,j_2,j_3)= \bigoplus_{k\in\cJ(j_1,j_2,j_3)}\cK_k \overline{AW}(j_1,j_2,j_3)\cK_k\ .
 \end{equation}
Then, confirming inequality \eqref{eq:min1} amounts to proving the following inequalities, for $k\in \cJ(j_1,j_2,j_3)$,
\begin{equation}
 \text{dim}(\cK_k \overline{AW}(j_1,j_2,j_3)\cK_k ) \leq  d_k^2  \ . \label{eq:min2}
\end{equation}
The algebras $\cK_k \overline{AW}(j_1,j_2,j_3)\cK_k$ are simpler to study than $\overline{AW}(j_1,j_2,j_3)$.
Roughly speaking, they correspond to replacing in the defining relations of $\overline{AW}(j_1,j_2,j_3)$ the central elements $\alpha_i$ (resp. $K$) by $\chi_{j_i}$ (resp. $\chi_k$).
One thus gets two annihilating polynomials for $A$ (similarly for $B$, $D$ and $D'$)
\begin{equation}
 \prod_{j\in\cJ(j_1,j_2)}(A-\chi_j)=0  \ , \qquad \prod_{m\in\cM(j_1,j_2,j_3)}(A-\chi_k+m)=0 \ ,
\end{equation}
which reduce to only one, {\it i.e.}
\begin{equation}
 \prod_{j\in\cJ^k(j_1,j_2,j_3)}(A-\chi_j)=0 \ , \label{eq:Aproj}
\end{equation}
with $\cJ^k(j_a,j_b,j_c)=\{j\in \cJ(j_a,j_b)\ | \ \chi_j\in \{\chi_k-m \ | \  m\in\cM(j_a,j_b,j_c) \}   \}$.

In fact, in the quotient of $\cC_{j_1,j_2,j_3}$ where $C_{123}=\chi_k$, the minimal polynomial of $C_{12}$ is 
\begin{equation}
\prod_{j\in\cS^k(j_1,j_2,j_3)}(C_{12}-\chi_j)=0\ ,
\end{equation}
where we recall that (see proof of Proposition \ref{pr:surj}) 
\begin{equation}
\cS^k(j_a,j_b,j_c) = \{\  j \in \cJ (j_a,  j_b )\ |\ k\in \cJ(j,j_c) \ \}\ .
\end{equation}
Similar results hold for $C_{23}$, $C_{13}^{(0)}$ and $C_{13}^{(1)}$. Let us emphasize that 
\begin{equation}
\cS^k(j_a,j_b,j_c) \subseteq \cJ^k(j_a,j_b,j_c) \ , \label{eq:ssu}
\end{equation}
and that the cardinality of $\cS^k(j_a,j_b,j_c)$ is equal to $d_k$ and does not depend on the ordering of $j_a,j_b,j_c$. 
This discussion suggests the definition of another quotient of $AW(3)$.
\begin{definition}
	The algebra $\overline{AW}^k(j_1,j_2,j_3)$, where $k\in\cJ(j_1,j_2,j_3)$, is the quotient of 
	the centrally extended Askey--Wilson algebra $AW(3)$ by $\alpha_i=\chi_{j_i}$ and the following relations
	\begin{align}
	& K=\chi_k\ ,\\
	&\prod_{j\in\cS^k(j_1,j_2,j_3)}(A-\chi_j)=0 \ , \prod_{j\in\cS^k(j_2,j_3,j_1)}(B-\chi_j)=0 \ ,  \label{eq:AWq11}\\
	&\prod_{j\in\cS^k(j_1,j_3,j_2)}\left(D-\chi_j \right)=0 \ , \prod_{j\in\cS^k(j_1,j_3,j_2)}\left(D'-\chi_j\right)=0 \ . \label{eq:AWq12}
	\end{align}
\end{definition}
\noindent Let us remark that the four annihilating polynomials \eqref{eq:AWq11}--\eqref{eq:AWq12} are of degree $d_k$.
These quotients lead to another conjecture.
\begin{conj} \label{conj2}
	The direct sum 
	\begin{equation}
	\widetilde{AW}(j_1,j_2,j_3)= \displaystyle
	\bigoplus_{k\in\cJ(j_1,j_2,j_3)}\overline{AW}^k(j_1,j_2,j_3)\ 
	\end{equation}
	is isomorphic to $\cC_{j_1,j_2,j_3}$.
\end{conj}
\noindent As for Conjecture \ref{conj1}, the proof of this conjecture reduces to showing that 
\begin{equation}
\text{dim}\left (\overline{AW}^k(j_1,j_2,j_3) \right)\leq d_k^2 \ . \label{eq:dim2}
\end{equation}
In view of \eqref{eq:ssu}, we see that Conjecture \ref{conj2} is true if Conjecture \ref{conj1} is.
Moreover, in this case, $\overline{AW}^k(j_1,j_2,j_3)$ is isomorphic to $\cK_k \overline{AW}(j_1,j_2,j_3)\cK_k$.

To conclude this section, let us emphasize that both conjectures presented above would provide an algebraic description of the centralizer $\cC_{j_1,j_2,j_3}$.
A strategy to prove these conjectures would be to establish inequalities \eqref{eq:dim2} and then to derive the isomorphism between $\overline{AW}^k(j_1,j_2,j_3)$ and $\cK_k \overline{AW}(j_1,j_2,j_3)\cK_k$.

\subsection{Invariance under permutations of $\{j_1,j_2,j_3\}$}\label{ssec:perm}

The algebras involved in Conjecture \ref{conj1} depend on the choice of three spins $j_1,j_2$ and $j_3$. We now show that it is sufficient to check the conjecture for only one ordering for the spins $j_1,j_2$ and $j_3$.   
\begin{prop}
	Let $j_1$, $j_2$ and $j_3$ be three positive half-integers or integers. If Conjecture 4.1 is true for the sequence of spins $\{j_1,j_2,j_3\}$, then it is also true for every permutation of $j_1,j_2,j_3$.
\end{prop}
\proof For any two representation maps $\pi_{j_1}$ and $\pi_{j_2}$ of $\Usl$, it is known that there exists an invertible matrix $P$ such that for all $x\in\Usl$ we have $(\kpd{\pi_{j_2}}{\pi_{j_1}})(\Delta(x))=P^{-1}(\kpd{\pi_{j_1}}{\pi_{j_2}})(\Delta(x))P$. Therefore, from the definition of the centralizer \eqref{eq:Cent} and the coassociativity of the comultiplication \eqref{eq:com2}, we deduce that for any permutation $\sigma$ of the symmetric group $\mathfrak{S}_3$, $\cC_{j_1,j_2,j_3}$ is isomorphic to $\cC_{j_{\sigma(1)},j_{\sigma(2)},j_{\sigma(3)}}$.

We must now show that the quotiented Askey--Wilson algebra $\overline{AW}(j_{\sigma(1)},j_{\sigma(2)},j_{\sigma(3)})$ is isomorphic to $\overline{AW}(j_1,j_2,j_3)$ for any permutation $\sigma\in\mathfrak{S}_3$. Since $\mathfrak{S}_3$ is generated by the transpositions $(1,2)$ and $(1,3)$, it suffices to prove the isomorphism for these two transformations.
  
The following maps are algebra isomorphisms
\begin{align}
	&\phi_1:\overline{AW}(j_3,j_2,j_1)\to\overline{AW}(j_1,j_2,j_3) \nonumber \\
	&\phi_1(\alpha_1)=\alpha_3, \ 
	 \phi_1(\alpha_2)=\alpha_2, \ 
	 \phi_1(\alpha_3)=\alpha_1, \ 
	 \phi_1(A)=B, \ 
	 \phi_1(B)=A, \ 
	 \phi_1(K)=K \ ,
\end{align}
\begin{align}
	&\phi_2:\overline{AW}(j_2,j_1,j_3)\to\overline{AW}(j_1,j_2,j_3) \nonumber \\
	&\phi_2(\alpha_1)=\alpha_2, \ 
	 \phi_2(\alpha_2)=\alpha_1, \ 
	 \phi_2(\alpha_3)=\alpha_3, \ 
	 \phi_2(A)=A, \ 
	 \phi_2(B)=D', \ 
	 \phi_2(K)=K \ .
\end{align}
To see the homomorphism for the defining relations of $AW(3)$, it is easier to work with the symmetric presentations \eqref{eq:AWcBD}--\eqref{eq:AWcBA} and \eqref{eq:AWcD'B}--\eqref{eq:AWcAD'}. By noticing that $\phi_1(D)=D'$ and $\phi_2(D)=B$, the homomorphism immediately follows. In order to preserve relation \eqref{eq:AWq0}, the central elements $\alpha_i$ have to be permuted in the same way as the spins $j_i$, which is the case for the two maps given above. For the quotiented relations \eqref{eq:AWq1}--\eqref{eq:AWq4}, the homomorphism is checked by observing that $\cJ(j_a,j_b)=\cJ(j_b,j_a)$, $\cM(j_a,j_b,j_c)=\cM(j_b,j_a,j_c)$ and that $\cJ(j_a,j_b,j_c)$ 
is invariant under any permutation of its entries. The R.H.S. of relation \eqref{eq:relOm} is invariant under any permutation of the elements $\alpha_i$. By using relations \eqref{eq:AWcBD}--\eqref{eq:AWcBA} and \eqref{eq:AWcD'B}--\eqref{eq:AWcAD'}, it is straightforward to show that $\phi_1(\Omega)=\Omega$ and $\phi_2(\Omega)=\Omega$, which proves the homomorphism for relation \eqref{eq:relOm}. Finally, since the maps $\phi_1$ and $\phi_2$ are surjective and invertible, they are bijective.
\endproof

\section{Finite irreducible representations of $\overline{AW}(j_1,j_2,j_3)$}\label{sec:repAW}
To support Conjecture \ref{conj1}, we want to show that the sum of the squares of the dimensions of all the finite irreducible representations of $\overline{AW}(j_1,j_2,j_3)$ is equal to the dimension of the centralizer. 
This result implies that Conjecture \ref{conj1} is true if and only if $\overline{AW}(j_1,j_2,j_3)$ is semisimple. Moreover, if $\overline{AW}(j_1,j_2,j_3)$ is not semisimple, the previous result proves that the missing relations (if there are any) in the kernel of $\phi$ are in a nilpotent radical of $\overline{AW}(j_1,j_2,j_3)$.

To identify all the finite irreducible representations of $\overline{AW}(j_1,j_2,j_3)$, we use the classification of the representations of the universal Askey--Wilson algebra given in \cite{H2} and look for the ones where the different relations of the quotient are satisfied. The universal Askey--Wilson algebra $\Delta_q$, introduced in \cite{Ter}, is generated by three elements $A,B,C$ and has three central elements $\alpha, \beta, \gamma$. There is a surjective algebra homomorphism from $\Delta_q$ to the quotient of $AW(3)$ by $\alpha_i=\chi_{j_i}$, with the following mappings 
\begin{gather}
A \mapsto A \ , \quad B \mapsto B \ , \quad C \mapsto D \ , \label{eq:homoABCAWuni} \\ 
\alpha \mapsto \chi_{j_1}\chi_{j_2}+\chi_{j_3}K \ , \quad \beta \mapsto \chi_{j_2}\chi_{j_3}+\chi_{j_1}K, \quad \gamma \mapsto \chi_{j_1}\chi_{j_3}+\chi_{j_2}K \ . \label{eq:homocentAWuni}
\end{gather}
We deduce that the quotient of $AW(3)$ by $\alpha_i=\chi_{j_i}$ is isomorphic to the quotient of $\Delta_q$ by the relations $(\alpha-\chi_{j_1}\chi_{j_2})/\chi_{j_3}=(\beta-\chi_{j_2}\chi_{j_3})/\chi_{j_1}=(\gamma-\chi_{j_1}\chi_{j_3})/\chi_{j_2}$. We can therefore use the representation theory of $\Delta_q$ in order to determine all the finite irreducible representations of $\overline{AW}(j_1,j_2,j_3)$.

The finite irreducible modules of the universal Askey--Wilson algebra $\Delta_q$ for $q$ not a root of unity are classified in \cite{H2}. They are given by the isomorphism classes of the $n+1$-dimensional modules $V_n(a,b,c)$ defined in \cite{H2}, for $n\geq0$, under certain conditions on $a,b,c$ (see Theorem 4.7 in \cite{H2}). In the representation $V_n(a,b,c)$, the central elements $\alpha,\beta$ and $\gamma$ of $\Delta_q$ take the following values  
\begin{align}
\alpha&=(q^{n+1}+q^{-n-1})(a+a^{-1})+(b+b^{-1})(c+c^{-1}) \ ,\label{eq:repal}\\
\beta&=(q^{n+1}+q^{-n-1})(b+b^{-1})+(c+c^{-1})(a+a^{-1}) \ ,\label{eq:repbe}\\
\gamma&=(q^{n+1}+q^{-n-1})(c+c^{-1})+(a+a^{-1})(b+b^{-1}) \ .\label{eq:repga}
\end{align}
The characteristic polynomials of $A,B,C \in \Delta_q$ in this representation are (see Lemma 4.3 of \cite{H2}) $K_a(X),K_b(X),K_c(X)$, with
\begin{equation}
K_x(X)=\prod_{i=0}^{n}(X-(q^{2i-n}x+q^{n-2i}x^{-1})) \ .\label{eq:reppolcar}
\end{equation}

The Casimir element of the algebra $\Delta_q$ is given by
\begin{equation}
qA\alpha +q^{-1} B \beta + q C\gamma -q^2 A^2 - q^{-2} B^2 -q^2 C^2 - q ABC \ . \label{eq:CasimirAWuni}
\end{equation}
It is straightforward to compute the value $\omega$ of this element in the representation $V_n(a,b,c)$ by using the representation matrices given in \cite{H2}. One gets  
\begin{align}
\omega&=(q^{n+1}+q^{-n-1})^2+(a+a^{-1})^2+(b+b^{-1})^2+(c+c^{-1})^2 \nonumber \\
&\quad +(q^{n+1}+q^{-n-1})(a+a^{-1})(b+b^{-1})(c+c^{-1})-(q+q^{-1})^2 \ . \label{eq:omAWuni}
\end{align}

We want to find all the irreducible representations of $\Delta_q$ that pass to the quotient $\overline{AW}(j_1,j_2,j_3)$. By comparing the annihilating polynomials of the elements $A,B,D \in \overline{AW}(j_1,j_2,j_3)$ given in \eqref{eq:AWq1}--\eqref{eq:AWq2} with the characteristic polynomials \eqref{eq:reppolcar}, we get the following restrictions for the inequivalent representations $V_n(a,b,c)$ to pass to the quotient :
\begin{align}
&0\leq n\leq \min\{j_1+j_2-|j_1-j_2|, j_2+j_3-|j_2-j_3|,j_1+j_3-|j_1-j_3|\} \ , \label{eq:repn}\\
&a=q^{2x+n+1} \ , \quad b=q^{2y+n+1} \ , \quad c=q^{2z+n+1} \ , \qquad \text{for } x,y,z \ \text{integers or half-integers,}   \label{eq:repabc}  \\ 
&|j_1-j_2| \leq x \leq j_1+j_2-n \ , \quad |j_2-j_3| \leq y \leq j_2+j_3-n \ , \quad |j_1-j_3| \leq z \leq j_1+j_3-n \ . \label{eq:repxyz} 
\end{align}

We recall that $K$ is a central element with the annihilating polynomial given in \eqref{eq:AWq1}. Therefore, $K$ has to be a constant equal to $\chi_\ell$ for some $\ell \in \cJ(j_1,j_2,j_3)$ in any irreducible representation of $\overline{AW}(j_1,j_2,j_3)$. We also recall that the Casimir element $\Omega$ of $AW(3)$ satisfies relation \eqref{eq:relOm} in the quotient $\overline{AW}(j_1,j_2,j_3)$. From this discussion and from the results \eqref{eq:homocentAWuni}--\eqref{eq:repga}, \eqref{eq:CasimirAWuni}, \eqref{eq:omAWuni} and \eqref{eq:repabc}, we deduce that the following equations must hold so that the representation $V_n(a,b,c)$ passes to the quotient $\overline{AW}(j_1,j_2,j_3)$ :
\begin{align}
&\chi_{\frac{n}{2}}\chi_{x+\frac{n}{2}}+\chi_{y+\frac{n}{2}}\chi_{z+\frac{n}{2}}
=\chi_{j_1}\chi_{j_2}+\chi_{j_3}\chi_{\ell}\ ,\label{eq:eqal}\\
&\chi_{\frac{n}{2}}\chi_{y+\frac{n}{2}}+\chi_{z+\frac{n}{2}}\chi_{x+\frac{n}{2}}
=\chi_{j_2}\chi_{j_3}+\chi_{j_1}\chi_{\ell}\ ,\label{eq:eqbe}\\
&\chi_{\frac{n}{2}}\chi_{z+\frac{n}{2}}+\chi_{x+\frac{n}{2}}\chi_{y+\frac{n}{2}}
=\chi_{j_1}\chi_{j_3}+\chi_{j_2}\chi_{\ell}\ ,\label{eq:eqga}\\
&\chi_{\frac{n}{2}}^2+\chi_{x+\frac{n}{2}}^2+\chi_{y+\frac{n}{2}}^2+\chi_{z+\frac{n}{2}}^2+\chi_{\frac{n}{2}}\chi_{x+\frac{n}{2}}\chi_{y+\frac{n}{2}}\chi_{z+\frac{n}{2}} =\chi_{j_1}^2+\chi_{j_2}^2+\chi_{j_3}^2+\chi_{\ell}^2+\chi_{j_1}\chi_{j_2}\chi_{j_3}\chi_{\ell} \ . \label{eq:eqom}
\end{align}

In the case of three identical spins $j_1=j_2=j_3=s$, we find by using mathematical software that there are only $192$ possible solutions for $x,y,z,n$ to the system of equations \eqref{eq:eqal}--\eqref{eq:eqom}. The only solutions respecting conditions \eqref{eq:repn}--\eqref{eq:repxyz} and corresponding to inequivalent representations $V_n(a,b,c)$ are
\begin{align}
&n=2\ell \ , &  &x=y=z=s-\ell \ ,  & &\text{if} \quad \ell \leq s \ , \label{eq:sol1} \\
&n = 3 s - \ell \ , & &x = y = z = \ell - s \ , & &\text{if} \quad \ell > s \ , \label{eq:sol2} \\
&n = s + \ell \ ,  & & x = y = 0, z = s - \ell \ ,  & & \text{if} \quad \ell < s \ , \label{eq:sol3}\\
&n = s - \ell - 1 \ ,  &  & x = y = 0 \ , z = s + l + 1 \ ,  & & \text{if} \quad \ell \leq s - 1 \label{eq:sol4} \ ,
\end{align}
and any permutation of $x,y,z$ in the previous equations is also a solution. 

Since $K=\chi_{\ell}$ for $\ell \in \cJ(j_1,j_2,j_3)$ in some irreducible representation, the annihilating polynomial of $K-A$ given in \eqref{eq:AWq3} implies that the annihilating polynomial of $A$ reduces to the relation \eqref{eq:Aproj} in this representation. If the set $\cJ^\ell(j_1,j_2,j_3)$ is equal to the set $\cS^\ell(j_1,j_2,j_3)$, then this reduced annihilating polynomial for $A$ leads to the constraint 
\begin{equation}
\text{max}( |j_1-j_2| , |j_3-\ell| )\leq x \leq \text{min}( j_1+j_2 , j_3+\ell )-n \ .
\end{equation}
Similar results hold for the annihilating polynomials of $B$ (resp. $D$) and the constraints on $y$ (resp. $z$). In the case of identical spins $j_1=j_2=j_3=s$, this implies that $s-\ell\leq x,y,z\leq s+\ell-n$ for $\ell\leq s$, and $\ell-s\leq x,y,z\leq 2s-n$ for $\ell>s$. The only solutions remaining are \eqref{eq:sol1} and \eqref{eq:sol2}. For $s$ half-integer, we do not find any cases where $\cJ^\ell(s,s,s)\neq\cS^\ell(s,s,s)$. For $s$ integer, as a consequence of the fact that $0\in\cM(s,s,s)$, we find $\cJ^\ell(s,s,s)\setminus\cS^\ell(s,s,s)=\{\ell\}$ if $\ell<s/2$, and otherwise the previous set is empty. We have verified numerically the sets $\cJ^\ell(s,s,s)\setminus\cS^\ell(s,s,s)$ given above for at least $s=\frac{1}{2},1,...,10$. In any case, the upper bound on the values of $x,y,z$ remains the same, and we still conclude that the only solutions are \eqref{eq:sol1} and \eqref{eq:sol2}. Therefore, the sum of the squares of the dimensions $n+1$ of all the irreducible modules of $\overline{AW}(s,s,s)$ is
\begin{equation}
	\sum_{\substack{\ell \in \cJ(s,s,s)\\
		\ell\leq s}} (2\ell+1)^2+ \sum_{\substack{\ell \in \cJ(s,s,s)\\
		\ell> s}} (3s-\ell+1)^2 = \frac{1}{2}(2s+1)((2s+1)^2+1) \ ,
\end{equation}
which is equal to the dimension of the centralizer $\dim(\cC_{s,s,s})=\displaystyle \sum_{\ell\in \cJ(s,s,s)} d_\ell^2$.
  
Let us remark that for $j_1=j_2=j_3=\frac{1}{2},1,...,\frac{17}{2}$, we used mathematical software to test all the possible integer values for $x,y,z,n$ such that the restrictions \eqref{eq:repn} and \eqref{eq:repxyz} and the three equations \eqref{eq:eqal}--\eqref{eq:eqga} are respected. The only solutions we found are those given in \eqref{eq:sol1}--\eqref{eq:sol4}. Hence, equation \eqref{eq:eqom} is perhaps not necessary if one wants to find all solutions for $x,y,z,n$ integers.

Let us also notice that in the general case where $j_1,j_2,j_3$ are any three fixed integers or half-integers, we find at least $192$ solutions to the system of equations \eqref{eq:eqal}--\eqref{eq:eqga} with $n$ integer and $x,y,z$ integers or half-integers. If these are the only such solutions, then it is possible to argue (in a similar manner as for the case of identical spins) that the sum of the squares of the dimensions of the irreducible representations that pass to the quotient $\overline{AW}(j_1,j_2,j_3)$ is also equal to the dimension of the centralizer $\cC_{j_1,j_2,j_3}$.

Finally, we notice that in the representations $V_n(a,b,c)$, the element $D'$ of $AW(3)$ has the same characteristic polynomial $K_c(X)$ as the element $D$. Therefore, the second relations in \eqref{eq:AWq2} and \eqref{eq:AWq4} do not provide any additional constraint on the values of $n,a,b,c$. 

\section{Quotient $\overline{AW}(\frac{1}{2},\frac{1}{2},\frac{1}{2})$ and Temperley--Lieb algebra}\label{sec:AW1TL}

In this section, we consider the case $j_1=j_2=j_3=\frac{1}{2}$ and show that the quotient $\overline{AW}(\frac{1}{2},\frac{1}{2},\frac{1}{2})$ is isomorphic to the centralizer $\cC_{\frac{1}{2},\frac{1}{2},\frac{1}{2}}$, which is known to be the Temperley--Lieb algebra. We give an explicit isomorphism between $\overline{AW}(\frac{1}{2},\frac{1}{2},\frac{1}{2})$ and the Temperley--Lieb algebra. 

\subsection{$\overline{AW}(\frac{1}{2},\frac{1}{2},\frac{1}{2})$ algebra}\label{ssec:AW1}
From the definitions \eqref{eq:J12}--\eqref{eq:J123} and \eqref{eq:M123}, we find the sets 
\begin{align}
	&\cJ\left(\frac{1}{2},\frac{1}{2}\right)=\{0,1\} \ , \quad \cJ\left(\frac{1}{2},\frac{1}{2},\frac{1}{2} \right)=\left\{\frac{1}{2},\frac{3}{2}\right\}, \label{eq:setJ1} \\ &\cM\left(\frac{1}{2},\frac{1}{2},\frac{1}{2} \right)=\{\chi_{\frac{3}{2}}-\chi_{1},\chi_{\frac{1}{2}}-\chi_{1},\chi_{\frac{1}{2}}-\chi_{0}\} \ . \label{eq:setM1}
\end{align}
The degeneracies are $d_{\frac{1}{2}}=2$ and $d_{\frac{3}{2}}=1$. We find from \eqref{eq:dimCent} that $\text{dim}\left(\cC_{\frac{1}{2},\frac{1}{2},\frac{1}{2}}\right)=5$. The central elements $\alpha_i$ of $\overline{AW}(\frac{1}{2},\frac{1}{2},\frac{1}{2})$ can all be replaced by the constant $\chi_{\frac{1}{2}}$. For computational convenience, we perform the transformation $X=(q-\qi)^2 \tilde{X}+q+\qi$ on the elements $X=A,B,D,D',K$, as in Remark \ref{rem:qlim}, and we define the shifted central element $\tG=\tK+[1/2]_q^2$. By using the sets given in \eqref{eq:setJ1} and \eqref{eq:setM1}, one finds that the defining relations of $\overline{AW}(\frac{1}{2},\frac{1}{2},\frac{1}{2})$ are
\begin{align}
&[\tB,[\tA,\tB]_q]_q=[2]_q\left(-\tB^2-\{\tA,\tB\}\right)+(q^2 + q^{-2})\tG\tB+[2]_q^2\tB \ ,\label{eq:AW1}\\
&[[\tA,\tB]_q,\tA]_q=[2]_q\left(-\tA^2-\{\tA,\tB\}\right)+(q^2 + q^{-2})\tG\tA+[2]_q^2\tA \ , \label{eq:AW2}\\
&\tA(\tA-[2]_q)=0 \ , \quad \tB(\tB-[2]_q)=0 \ , \quad (\tG-1)(\tG-[2]_q^2)=0 \ , \label{eq:AWABK}\\
&\tD(\tD-[2]_q)=0 \ , \quad \tD'(\tD'-[2]_q)=0 \ ,\label{eq:AWD}\\
&(\tG-\tA+[2]_q-[2]_q^2)(\tG-\tA+[2]_q-1)(\tG-\tA-1)=0 \ , \label{eq:AWKA}\\
&(\tG-\tB+[2]_q-[2]_q^2)(\tG-\tB+[2]_q-1)(\tG-\tB-1)=0 \ , \label{eq:AWKB}\\
&(\tG-\tD+[2]_q-[2]_q^2)(\tG-\tD+[2]_q-1)(\tG-\tD-1)=0 \ , \label{eq:AWKD}\\
&(\tG-\tD'+[2]_q-[2]_q^2)(\tG-\tD'+[2]_q-1)(\tG-\tD'-1)=0 \ , \label{eq:AWKD'}\\
&(q^2 + q^{-2})(q-\qi)(q\tA + \qi \tB + q\tD)\tG - [2]_q(([2]_q-q^3)(\tA + \tD) + q^{-3}\tB) \label{eq:AWOm} \\
&
- (q-\qi)(q^2\tA^2 + q^{-2}\tB^2 + q^2\tD^2)
- q[2]_q(q-\qi)(\tA\tB + \tA\tD + \tB\tD)
- q(q-\qi)^3\tA\tB\tD
 \nonumber \\
&=(q-\qi)\tG^2
+ (q^5 - q^{-5} - q^2 [2]_q)\tG
- q^{-1} [2]_q^2 \ , \nonumber
\end{align}
where
\begin{align}
\tD=[2]_q+(q^2+q^{-2})\tG-\tA-\tB-\frac{q-\qi}{q+\qi}[\tA,\tB]_q \ , \label{eq:relD}\\
\tD'=[2]_q+(q^2+q^{-2})\tG-\tA-\tB-\frac{q-\qi}{q+\qi}[\tB,\tA]_q \ .	
\end{align}
We want to show that $\overline{AW}(\frac{1}{2},\frac{1}{2},\frac{1}{2})$ is isomorphic to the centralizer $\cC_{\frac{1}{2},\frac{1}{2},\frac{1}{2}}$.
\begin{prop}
	The relations defining the quotient $\overline{AW}(\frac{1}{2},\frac{1}{2},\frac{1}{2})$ can be given as follows
	\begin{align}
	&\tA^2=[2]_q \tA \ , \quad \tB^2=[2]_q \tB \ , \label{eq:AB} \\ 
	&\tA\tB\tA=[2]_q \{\tA,\tB\}-[3]_q\tA-[2]_q^2\tB+[2]_q[3]_q \ , \label{eq:ABA} \\
	&\tB\tA\tB=[2]_q \{\tA,\tB\}-[3]_q\tB-[2]_q^2\tA+[2]_q[3]_q \ . \label{eq:BAB}
	\end{align}
\end{prop}
\proof
The two first relations in \eqref{eq:AWABK} directly lead to \eqref{eq:AB}. The third relation in \eqref{eq:AWABK} implies  
\begin{equation}
\tG^2=([2]_q^2+1)\tG-[2]_q^2 \ .\label{eq:AWG}
\end{equation}
Developing \eqref{eq:AW1} and \eqref{eq:AW2} and using \eqref{eq:AB}, one gets 
\begin{equation}
\tB\tA\tB=\tG\tB \ , \quad \tA\tB\tA=\tG\tA \ . \label{eq:GAGB}
\end{equation}
Expanding \eqref{eq:AWKA} and \eqref{eq:AWKB} and simplifying with the help of \eqref{eq:AB} and \eqref{eq:AWG}, one gets
\begin{equation}
	\tG\tA =[2]_q \tG +\tA-[2]_q \ , \quad \tG\tB =[2]_q \tG +\tB-[2]_q \ , \label{eq:GAGB2}
\end{equation}
which implies
\begin{align}
\tG\tA\tB =[2]_q^2 \tG+\tA\tB-[2]_q^2 \ , \quad
\tG\tB\tA =[2]_q^2 \tG+\tB\tA-[2]_q^2  \ . \label{eq:GABGBA}
\end{align}
Equations \eqref{eq:AWD} and \eqref{eq:AWKD}--\eqref{eq:AWOm} can be simplified using the previous relations, and they lead to
\begin{equation}
	\tG=-[2]_q(\tA+\tB)+\{\tA,\tB\}+[2]_q^2 \ . \label{eq:G2}
\end{equation}
Substituting \eqref{eq:G2} in \eqref{eq:GAGB2} and \eqref{eq:GABGBA}, one finds
\begin{align}
&\tG\tA=[2]_q \{\tA,\tB\}-[3]_q\tA-[2]_q^2\tB+[2]_q[3]_q \ , \label{eq:GA} \\
&\tG\tB=[2]_q \{\tA,\tB\}-[3]_q\tB-[2]_q^2\tA+[2]_q[3]_q \ , \label{eq:GB}\\ 
&\tG\tA\tB=[2]_q^2\tB\tA+([2]_q^2+1)\tA\tB-[2]_q^3(\tA+\tB)+[2]_q^2[3]_q \ , \label{eq:GAB} \\
&\tG\tB\tA=[2]_q^2\tA\tB+([2]_q^2+1)\tB\tA-[2]_q^3(\tA+\tB)+[2]_q^2[3]_q \ . \label{eq:GBA}
\end{align}
Equations \eqref{eq:GAGB} and \eqref{eq:GA}--\eqref{eq:GB} imply the relations \eqref{eq:ABA}--\eqref{eq:BAB} of the proposition. 

It remains to show that the generator $\tG$ can be suppressed from the presentation, or in other words that \eqref{eq:AWG} and \eqref{eq:GA}--\eqref{eq:GBA} are implied from the relations of the proposition. Suppose that relations \eqref{eq:AB}--\eqref{eq:BAB} are true and let $\tG=-[2]_q(\tA+\tB)+\{\tA,\tB\}+[2]_q^2$. Multiplying the expression of $\tG$ on the left and on the right by $\tA$ and $\tB$, one finds
\begin{equation}
	\tG\tA=\tA\tG=\tA\tB\tA \ , \quad \tG\tB=\tB\tG=\tB\tA\tB \ .
\end{equation}
Using \eqref{eq:ABA} and \eqref{eq:BAB}, equations \eqref{eq:GA} and \eqref{eq:GB} are recovered. Multiplying \eqref{eq:ABA} on the right by $\tB$ and \eqref{eq:BAB} on the right by $\tA$, one finds
\begin{align}
 & \tG\tA\tB=\tA\tB\tA\tB=\tA\tB-[2]_q \tB +[2]_q \tB\tA\tB \ , \label{eq:ABAB}\\
 & \tG\tB\tA=\tB\tA\tB\tA=\tB\tA-[2]_q \tA +[2]_q \tA\tB\tA \ , \label{eq:BABA}
\end{align}
from which one easily recovers \eqref{eq:GAB} and \eqref{eq:GBA}. Finally, it is straightforward to arrive at 
\begin{align}
\tG^2=-[2]_q^3(\tA+\tB)+[2]_q^2\{\tA,\tB\}-[2]_q(\tA\tB\tA+\tB\tA\tB)+\tA\tB\tA\tB+\tB\tA\tB\tA+[2]_q^4
\end{align}
and to use the results \eqref{eq:ABAB} and \eqref{eq:BABA} to recover \eqref{eq:AWG}.
\endproof

\begin{thm}
	Conjecture \ref{conj1} is verified for $j_1=j_2=j_3=\frac{1}{2}$.
\end{thm}
\proof
We already know from proposition \ref{pr:surj} that the map $\overline{\phi}$ is surjective. From the previous proposition, it is easy to show that $\{1,\tA,\tB,\tA\tB,\tB\tA\}$ is a linearly generating set of $\overline{AW}(\frac{1}{2},\frac{1}{2},\frac{1}{2})$. Since $\text{dim}\left(\cC_{\frac{1}{2},\frac{1}{2},\frac{1}{2}}\right)=5$, this shows the injectivity of the map $\overline{\phi}$.
\endproof

\subsection{Connection with the Temperley--Lieb algebra}\label{ssec:TL}

It is known that the Temperley--Lieb algebra is isomorphic to the centralizer of the diagonal embedding of $\Usl$ in the tensor product of three fundamental representations \cite{Jimbo}. 
Hence, from the results of the previous subsection, the quotiented Askey--Wilson algebra $\overline{AW}(\frac{1}{2},\frac{1}{2},\frac{1}{2})$ is isomorphic to the Temperley--Lieb algebra.

\begin{definition}\cite{TL} 
	The Temperley--Lieb algebra $TL_3(q)$ is generated by $\sigma_1$ and $\sigma_2$ with the following defining relations
	\begin{align}
	&\sigma_1^2=(q+\qi)\sigma_1 \ , \quad \sigma_2^2=(q+\qi)\sigma_2 \ , \label{eq:TL1} \\
	&\sigma_1\sigma_2\sigma_1=\sigma_1 \ , \quad \sigma_2\sigma_1\sigma_2=\sigma_2 \ . \label{eq:TL2} 
	\end{align}
\end{definition}

\begin{thm}
	The quotiented Askey--Wilson algebra $\overline{AW}(\frac{1}{2},\frac{1}{2},\frac{1}{2})$ is isomorphic to the Temperley--Lieb algebra $TL_3(q)$. This isomorphism is given explicitly by
	\begin{align}
	\overline{AW}\left(\frac{1}{2},\frac{1}{2},\frac{1}{2}\right) &\to TL_3(q) \nonumber \\
	\tA&\mapsto(q+\qi)-\sigma_1 \ , \label{eq:isoTL1} \\
	\tB&\mapsto(q+\qi)-\sigma_2 \ . \label{eq:isoTL2} 
	\end{align}
\end{thm}
\proof
It is straightforward to show that the defining relations \eqref{eq:TL1} and \eqref{eq:TL2} of $TL_3(q)$ are equivalent to the relations \eqref{eq:AB}--\eqref{eq:BAB} of $\overline{AW}(\frac{1}{2},\frac{1}{2},\frac{1}{2})$.
\endproof
\section{Quotient $\overline{AW}(1,1,1)$ and Birman--Murakami--Wenzl algebra}\label{sec:AW2BMW}
In this section, we choose $j_1=j_2=j_3=1$ and prove that the quotient $\overline{AW}(1,1,1)$ is isomorphic to the centralizer $\cC_{1,1,1}$. In this case, $\cC_{1,1,1}$ is known to be connected to the Birman--Murakami--Wenzl algebra. We give an explicit isomorphism between $\overline{AW}(1,1,1)$ and a specialization of the BMW algebra.
\subsection{$\overline{AW}(1,1,1)$ algebra}\label{ssec:AW2}	
We have the following sets
\begin{gather}
\cJ(1,1)=\{0,1,2\} \ , \quad \cJ(1,1,1)=\{0,1,2,3\} \ , \\ \cM(1,1,1)=\{\chi_1-\chi_2,\chi_0-\chi_1,0,\chi_1-\chi_0,\chi_2-\chi_1,\chi_3-\chi_2\} \ .
\end{gather}
The degeneracies are $d_0=d_3=1$, $d_1=3$ and $d_2=2$, and the dimension of the centralizer is $\dim(C_{1,1,1})=15$. The central elements $\alpha_i$ of $\overline{AW}(1,1,1)$ can all be replaced by the constant $\chi_1$. For computational convenience again, we perform the transformation $X=(q-\qi)^2 \tilde{X}+q+\qi$ on the elements $X=A,B,D,D',K$ (see Remark \ref{rem:qlim}). We recall that the eigenvalues $\chi_j$ are transformed to $\tchi_j=[j]_q[j+1]_q$, and we define the constants $m_1=\tchi_1-\tchi_0$, $m_2=\tchi_2-\tchi_1$ and $m_3=\tchi_3-\tchi_2$. The defining relations \eqref{eq:AWq1}--\eqref{eq:AWq4} of $\overline{AW}(1,1,1)$ are written as
\begin{align}
&(q^2+q^{-2})\tB\tA\tB=\tA\tB^2+\tB^2\tA-[2]_q \tB^2-[2]_q\{\tA,\tB\}+[2]_q([2]_q^2-3)\tK\tB+[2]_q^2[3]_q\tB\label{eq:2AW1} \ ,\\
&(q^2+q^{-2})\tA\tB\tA=\tB\tA^2+\tA^2\tB-[2]_q \tA^2-[2]_q\{\tA,\tB\}+[2]_q([2]_q^2-3)\tK\tA+[2]_q^2[3]_q\tA\label{eq:2AW2} \ , \\
&\tA(\tA-\tchi_1)(\tA-\tchi_2)=0
 \ , \quad \tB(\tB-\tchi_1)(\tB-\tchi_2)=0
 \ , \label{eq:2AWAB} \\ 
&\tK(\tK-\tchi_1)(\tK-\tchi_2)(\tK-\tchi_3)=0 \ , \label{eq:2AWK} \\ 
&\tD(\tD-\tchi_1)(\tD-\tchi_2)=0 \ , \quad \tD'(\tD'-\tchi_1)(\tD'-\tchi_2)=0 \ , \label{eq:2AWD} \\ 
&(\tK-\tA+m_2)(\tK-\tA+m_1)(\tK-\tA)(\tK-\tA-m_1)(\tK-\tA-m_2)(\tK-\tA-m_3)=0 \ , \label{eq:2AWKA}\\
&(\tK-\tB+m_2)(\tK-\tB+m_1)(\tK-\tB)(\tK-\tB-m_1)(\tK-\tB-m_2)(\tK-\tB-m_3)=0 \ , \label{eq:2AWKB}\\
&(\tK-\tD+m_2)(\tK-\tD+m_1)(\tK-\tD)(\tK-\tD-m_1)(\tK-\tD-m_2)(\tK-\tD-m_3)=0 \ , \label{eq:2AWKD} \\
&(\tK-\tD'+m_2)(\tK-\tD'+m_1)(\tK-\tD')(\tK-\tD'-m_1)(\tK-\tD'-m_2)(\tK-\tD'-m_3)=0 \ , \label{eq:2AWKD'}
\end{align}
where
\begin{align}
\tD&=[2]_q[3]_q+([2]_q^2-3)\tK-\frac{(q-\qi)}{[2]_q}[\tA,\tB]_q-\tA-\tB \ ,\\
\tD'&=[2]_q[3]_q+([2]_q^2-3)\tK-\frac{(q-\qi)}{[2]_q}[\tB,\tA]_q-\tA-\tB \ .	
\end{align}

\begin{thm}
	Conjecture 4.1 is verified for $j_1=j_2=j_3=1$.
\end{thm}
\proof
We already know from proposition \ref{pr:surj} that the map $\overline{\phi}$ is surjective. We only need to prove that it is injective in this case.

We define the set
\begin{equation}
\mathcal{S}=\{1,\tA,\tB,\tA^2,\tA\tB,\tB\tA,\tB^2,\tA^2\tB,\tA\tB\tA,\tA\tB^2,\tB\tA^2,\tB\tA\tB,\tA^2\tB^2,\tA\tB\tA\tB,\tB\tA\tB\tA\} \ .
\end{equation}
Using relations \eqref{eq:2AW1}--\eqref{eq:2AWK}, it can be shown that $\mathcal{S}_r=\mathcal{S}\cup \tK\mathcal{S}\cup \tK^2\mathcal{S}\cup \tK^3\mathcal{S}$ is a linearly generating set for $\overline{AW}(1,1,1)$. We can construct the 60 by 60 matrices $\tA_r$, $\tB_r$ and $\tK_r$ corresponding to the regular actions of $\tA$, $\tB$ and $\tK$ on the set $\mathcal{S}_r$. Knowing that $\tA_r$, $\tB_r$ and $\tK_r$ have to satisfy \eqref{eq:2AW1}--\eqref{eq:2AWK} and the first relation of \eqref{eq:2AWD}, we find 32 independant relations between the elements of $\mathcal{S}_r$ and we can reduce the generating set to
\begin{equation}
\mathcal{S}'_r=\mathcal{S}\cup \tK\{1,\tB,\tA^2,\tA\tB,\tB^2,\tA^2\tB,\tA\tB\tA,\tA\tB^2,\tB\tA\tB,\tA^2\tB^2,\tA\tB\tA\tB\}\cup \tK^2\{\tA^2,\tA^2\tB\} \ .
\end{equation}
We repeat the procedure and construct 28 by 28 matrices corresponding to the regular actions on $\mathcal{S}'_r$. Only using again \eqref{eq:2AW1}--\eqref{eq:2AWK} and the first of \eqref{eq:2AWD}, we can reduce the generating set to  
\begin{equation}
\mathcal{S}''_r=\mathcal{S}\cup\{\tK,\tK\tB,\tK\tA^2,\tK\tA\tB,\tK\tB^2,\tK\tA^2\tB,\tK\tA\tB\tA,\tK\tA\tB^2,\tK\tB\tA\tB\} \ .
\end{equation}
We repeat and construct 24 by 24 matrices. At this point, relations \eqref{eq:2AW1}--\eqref{eq:2AWD} are already satisfied. We must use relations \eqref{eq:2AWKA}--\eqref{eq:2AWKD} to reduce the generating set to
\begin{equation}
\mathcal{S}'''_r=\mathcal{S}\cup\{\tK,\tK\tA^2,\tK\tB^2\} \ .
\end{equation}
We repeat one last time by constructing 18 by 18 matrices and we use \eqref{eq:2AW1}--\eqref{eq:2AWKD'} to find 3 independant relations which allow to reduce the generating set to $\mathcal{S}$. It can also be verified that the matrices of the regular action satisfy the defining relation \eqref{eq:relOm} involving the Casimir element $\Omega$. We made the previous computations by using a formal mathematical software.

From these results, we have that $\mathcal{S}$ is a linearly generating set for $\overline{AW}(1,1,1)$ with 15 elements. Since $\text{dim}\left(\cC_{1,1,1}\right)=15$, we conclude that $\overline{\phi}$ is injective.
\endproof

\subsection{Connection with the Birman--Murakami--Wenzl algebra}\label{ssec:BMW}

It is known \cite{LZ} that the Birman--Murakami--Wenzl algebra is isomorphic to the centralizer of the diagonal embedding of $\Usl$ in the tensor product of three spin-$1$ representations. 
Hence, from the previous theorem, the quotiented Askey--Wilson algebra $\overline{AW}(1,1,1)$ is isomorphic to the BMW algebra.

\begin{definition} \cite{IMO}
	The Birman--Murakami--Wenzl algebra $BMW_3(Q,\mu)$ is generated by invertible elements $s_1$ and $s_2$ with the following defining relations
	\begin{align}
	&s_1s_2s_1=s_2s_1s_2 \ , \label{eq:BMWdef1}\\
	&e_1s_1=s_1e_1=\mu^{-1} e_1 \ , \quad e_2s_2=s_2e_2=\mu^{-1} e_2 \ , \label{eq:BMWdef2} \\
	&e_1s_2^{\epsilon}e_1=\mu^{\epsilon}e_1 \ , \quad e_2s_1^{\epsilon}e_2=\mu^{\epsilon}e_2 \ , \quad \epsilon=\pm 1 \ , \label{eq:BMWdef3}\\
	&e_i=1-\frac{s_i-s_i^{-1}}{Q-Q^{-1}} \ , \quad i=1,2 \ . \label{eq:BMWdef4}
	\end{align}
\end{definition}

\begin{thm}
	The quotiented Askey--Wilson algebra $\overline{AW}(1,1,1)$ is isomorphic to the Birman--Murakami--Wenzl algebra $BMW_3(q^2,q^4)$. This isomorphism is given explicitly by
	\begin{align}
	\overline{AW}(1,1,1) &\to BMW_3(q^2,q^4) \nonumber \\
	\tA&\mapsto(q+\qi)(s_1-q^{-2} e_1)+(q+\qi)^2\qi \ , \label{eq:BMWiso1} \\
	\tB&\mapsto(q+\qi)(s_2-q^{-2} e_2)+(q+\qi)^2\qi \ . \label{eq:BMWiso2}
	\end{align}
\end{thm}
\proof
The algebras $\overline{AW}(1,1,1)$ and $BMW_3(q^2,q^4)$ are both isomorphic to $\cC_{1,1,1}$, hence they are isomorphic to each other.
It can be verified that the image of $\tA$ (resp. $\tB$) in $\End(\M{1}^{\otimes 3})$ is equal to the image of the R.H.S. of \eqref{eq:BMWiso1} (resp.  \eqref{eq:BMWiso2}), which justifies the explicit mapping. The inverse map is given by 
\begin{align}
s_1&\mapsto q^{-2}(q+\qi)^{-2}\tA^2-q^{-2}(q+\qi)^{-1}(2+q^{-2})\tA+q^{-4} \ , \label{eq:BMWiso3} \\
s_2&\mapsto q^{-2}(q+\qi)^{-2}\tB^2-q^{-2}(q+\qi)^{-1}(2+q^{-2})\tB+q^{-4} \ . \label{eq:BMWiso4}
\end{align}
\endproof

\section{Quotient $\overline{AW}(\frac{3}{2},\frac{3}{2},\frac{3}{2})$}\label{sec:AW3}
In this section, we take $j_1=j_2=j_3=\frac{3}{2}$ and show that the $\overline{AW}(\frac{3}{2},\frac{3}{2},\frac{3}{2})$ algebra is isomorphic to the centralizer $\cC_{\frac{3}{2},\frac{3}{2},\frac{3}{2}}$.

From the decomposition rules of the tensor product, we find the sets
\begin{equation}
\cJ\left(\frac{3}{2},\frac{3}{2}\right)=\{0,1,2,3\} \ , \quad \cJ\left(\frac{3}{2},\frac{3}{2},\frac{3}{2} \right)=\left\{\frac{1}{2},\frac{3}{2},\frac{5}{2},\frac{7}{2},\frac{9}{2}\right\} \ ,
\end{equation}
\begin{equation}
\begin{split}
\cM\left(\frac{3}{2},\frac{3}{2},\frac{3}{2} \right)=\{&\chi_{\frac{9}{2}}-\chi_{3},\chi_{\frac{7}{2}}-\chi_{3},\chi_{\frac{7}{2}}-\chi_{2},\chi_{\frac{5}{2}}-\chi_{3},\chi_{\frac{5}{2}}-\chi_{2},\chi_{\frac{5}{2}}-\chi_{1},\\
&\chi_{\frac{3}{2}}-\chi_{3},\chi_{\frac{3}{2}}-\chi_{2},\chi_{\frac{3}{2}}-\chi_{1},\chi_{\frac{3}{2}}-\chi_{0},\chi_{\frac{1}{2}}-\chi_{2},\chi_{\frac{1}{2}}-\chi_{1}\} \ .
\end{split}
\end{equation}
The degeneracies are $d_{\frac{9}{2}}=1$, $d_{\frac{7}{2}}=d_{\frac{1}{2}}=2$, $d_{\frac{5}{2}}=3$ and $d_{\frac{3}{2}}=4$, and the dimension of the centralizer is $\dim\left(\cC_{\frac{3}{2},\frac{3}{2},\frac{3}{2}}\right)=34$. The central elements $\alpha_i$ of $\overline{AW}(\frac{3}{2},\frac{3}{2},\frac{3}{2})$ are each equal to the constant $\chi_{\frac{3}{2}}$. 

In order to prove the injectivity of the map $\overline{\phi}$ in this case, we will use the strategy described in Subsection 4.3 and show that
\begin{equation}
\text{dim}\left(\cK_k \overline{AW}\left(\frac{3}{2},\frac{3}{2},\frac{3}{2}\right)\cK_k \right) \leq  d_k^2  \qquad \forall k\in \cJ\left(\frac{3}{2},\frac{3}{2},\frac{3}{2}\right)\ . \label{eq:dimproject}
\end{equation}
We recall that for each $k\in \cJ\left(\frac{3}{2},\frac{3}{2},\frac{3}{2}\right)$, the central element $K$ is replaced by the constant $\chi_k$ in the algebra $\cK_k \overline{AW}\left(\frac{3}{2},\frac{3}{2},\frac{3}{2}\right)\cK_k$, and the annihilating polynomials for $A$ (similarly for $B$, $D$ and $D'$) reduce to
\begin{equation}
\prod_{j\in\cJ^k\left(\frac{3}{2},\frac{3}{2},\frac{3}{2}\right)}(A-\chi_j)=0 \ ,
\end{equation}
where $\cJ^k\left(\frac{3}{2},\frac{3}{2},\frac{3}{2}\right)=\{j\in \cJ\left(\frac{3}{2},\frac{3}{2}\right)\ | \ \chi_j\in \{\chi_k-m \ | \  m\in\cM\left(\frac{3}{2},\frac{3}{2},\frac{3}{2}\right) \}   \}$. Once again, we perform the transformation $X=(q-\qi)^2 \tilde{X}+q+\qi$ on the elements $X=A,B,D,D',K$ (see Remark \ref{rem:qlim}). Therefore, one finds that the following relations hold in the algebras $\cK_k \overline{AW}\left(\frac{3}{2},\frac{3}{2},\frac{3}{2}\right)\cK_k$ :
\begin{itemize}
	\item $k=\frac{9}{2}$ 
	\begin{equation}
	\tA=\tB=\chi_{3} \ .
	\end{equation}
	\item $k=\frac{7}{2}$
	\begin{align}
	&\tB\tA\tB=[2]_q^3\{\tA,\tB\}+[3]_q^2(-2[2]_q^2\tA+(q^4+q^{-4}-1)\tB+[2]_q^3) \ , \\
	&\tA\tB\tA=[2]_q^3\{\tA,\tB\}+[3]_q^2(-2[2]_q^2\tB+(q^4+q^{-4}-1)\tA+[2]_q^3) \ , \\
	&(\tA-\tchi_2)(\tA-\tchi_3)=0
	\ , \quad (\tB-\tchi_2)(\tB-\tchi_3)=0
	\ , \\
	&(\tD-\tchi_2)(\tD-\tchi_3)=0 \ , \quad (\tD'-\tchi_2)(\tD'-\tchi_3)=0 \ . \label{eq:relDproj1}
	\end{align}
	\item $k=\frac{5}{2}$
	\begin{align}
	&(q^2+q^{-2})\tB\tA\tB=\tA\tB^2+\tB^2\tA-[2]_q \tB^2-[2]_q\{\tA,\tB\}+(2[2]_q\tchi_{\frac{3}{2}}+(q^4+q^{-4})(\tchi_{\frac{3}{2}}+\tchi_{\frac{5}{2}}))\tB \ , \\
	&(q^2+q^{-2})\tA\tB\tA=\tB\tA^2+\tA^2\tB-[2]_q \tA^2-[2]_q\{\tA,\tB\}+(2[2]_q\tchi_{\frac{3}{2}}+(q^4+q^{-4})(\tchi_{\frac{3}{2}}+\tchi_{\frac{5}{2}}))\tA \ , \\
	&(\tA-\tchi_1)(\tA-\tchi_2)(\tA-\tchi_3)=0
	\ , \quad (\tB-\tchi_1)(\tB-\tchi_2)(\tB-\tchi_3)=0
	\ , \\
	&(\tD-\tchi_1)(\tD-\tchi_2)(\tD-\tchi_3)=0 \ , \quad (\tD'-\tchi_1)(\tD'-\tchi_2)(\tD'-\tchi_3)=0 \ .
	\end{align}
	\item $k=\frac{3}{2}$
	\begin{align}
	&(q^2+q^{-2})\tB\tA\tB=\tA\tB^2+\tB^2\tA-[2]_q \tB^2-[2]_q\{\tA,\tB\}+2([2]_q+q^4+q^{-4})\tchi_{\frac{3}{2}}\tB \ , \\
	&(q^2+q^{-2})\tA\tB\tA=\tB\tA^2+\tA^2\tB-[2]_q \tA^2-[2]_q\{\tA,\tB\}+2([2]_q+q^4+q^{-4})\tchi_{\frac{3}{2}}\tA \ , \\
	&\tA(\tA-\tchi_1)(\tA-\tchi_2)(\tA-\tchi_3)=0
	\ , \quad \tB(\tB-\tchi_1)(\tB-\tchi_2)(\tB-\tchi_3)=0
	\ , \\
	&\tD(\tD-\tchi_1)(\tD-\tchi_2)(\tD-\tchi_3)=0 \ , \quad \tD'(\tD'-\tchi_1)(\tD'-\tchi_2)(\tD'-\tchi_3)=0 \ .
	\end{align}
	\item $k=\frac{1}{2}$
	\begin{align}
	&(q^2+q^{-2})\tB\tA\tB=[3]_q([2]_q\{\tA,\tB\}-2[2]_q^2\tA+[2]_q^3)+[2]_q^2(q^4+q^{-4})\tB \ , \\
	&(q^2+q^{-2})\tA\tB\tA=[3]_q([2]_q\{\tA,\tB\}-2[2]_q^2\tB+[2]_q^3)+[2]_q^2(q^4+q^{-4})\tA \ , \\
	&(\tA-\tchi_1)(\tA-\tchi_2)=0
	\ , \quad (\tB-\tchi_1)(\tB-\tchi_2)=0
	\ , \\
	&(\tD-\tchi_1)(\tD-\tchi_2)=0 \ , \quad (\tD'-\tchi_1)(\tD'-\tchi_2)=0 \ . \label{eq:relDproj2}
	\end{align}
\end{itemize}
For each value of $k$, the elements $\tD$ and $\tD'$ are given by
\begin{align}
\tD&=\frac{(q^4+q^{-4})}{[2]_q}(\tchi_{\frac{3}{2}}+\tchi_k)+2\tchi_{\frac{3}{2}}-\frac{(q-\qi)}{[2]_q}[\tA,\tB]_q-\tA-\tB \ , \label{eq:Dproject}\\
\tD'&=\frac{(q^4+q^{-4})}{[2]_q}(\tchi_{\frac{3}{2}}+\tchi_k)+2\tchi_{\frac{3}{2}}-\frac{(q-\qi)}{[2]_q}[\tB,\tA]_q-\tA-\tB \ .	
\end{align}

\begin{thm}
	Conjecture 4.1 is verified for $j_1=j_2=j_3=\frac{3}{2}$.
\end{thm}
\proof
Since we already know that the map $\overline{\phi}$ is surjective, we only need to prove \eqref{eq:dimproject}. For the case $k=\frac{9}{2}$, all the elements are constants and $\dim(\cK_\frac{9}{2} \overline{AW}\left(\frac{3}{2},\frac{3}{2},\frac{3}{2}\right)\cK_\frac{9}{2})=1$. For the case $k=\frac{7}{2}$ (resp. $k=\frac{1}{2}$), one uses \eqref{eq:Dproject} in the first relation of \eqref{eq:relDproj1} (resp. \eqref{eq:relDproj2}) to find (after some simplifications using the defining relations ) $\tB\tA=-\tA\tB+x_1(\tA+\tB)+x_2$, for some constants $x_1$ and $x_2$ that can be computed. Therefore, in both cases we see that a linearly generating set is given by $\{1,\tA,\tB,\tA\tB\}$. For the case $k=\frac{5}{2}$, we used formal mathematical software to show that $\{1,\tA,\tB,\tA^2,\tA\tB,\tB\tA,\tB^2,\tA\tB^2,\tB\tA^2\}$ is a generating set. Similarly, for the case $k=\frac{3}{2}$, a generating set is given by $\{1,\tA,\tB,\tA^2,\tA\tB,\tB\tA,\tB^2,\tA^3,\tA\tB\tA,\tA\tB^2,\tB\tA^2,\tB\tA\tB,\tB^3,\tA^3\tB,\tA^2\tB^2,\tA\tB^3\}$. From these results and the degeneracies $d_k$ given at the beginning of the section, we see that \eqref{eq:dimproject} holds, which concludes the proof.
\endproof

Let us notice that the defining relation \eqref{eq:relOm} of $\overline{AW}(j_1,j_2,j_3)$ which involves the Casimir element $\Omega$ of $AW(3)$ has not been called upon in the previous proof. It is straightforward to verify that this relation is satisfied in each of the algebras $\cK_k \overline{AW}\left(\frac{3}{2},\frac{3}{2},\frac{3}{2}\right)\cK_k$ by using the relations given above.

\section{Quotient $\overline{AW}(j,\frac{1}{2},\frac{1}{2})$ and one-boundary Temperley--Lieb algebra}\label{sec:AWj1bTL}
In this section, we consider the case $j_1=j$, for $j=1,\frac{3}{2},2,...$, and $j_2=j_3=\frac{1}{2}$. We show that the algebra $\overline{AW}(j,\frac{1}{2},\frac{1}{2})$ is isomorphic to the centralizer $\cC_{j,\frac{1}{2},\frac{1}{2}}$. We also find an explicit isomorphism between this quotient of the Askey--Wilson algebra and a specialization of the one-boundary Temperley--Lieb algebra.
\subsection{$\overline{AW}(j,\frac{1}{2},\frac{1}{2})$ algebra}\label{ssec:AWj}

From the tensor decomposition rules, we find the sets
\begin{align*}
&\cJ\left(j,\frac{1}{2}\right)=\left\{j-\frac{1}{2},j+\frac{1}{2}\right\} \ , \quad \cJ\left(\frac{1}{2},\frac{1}{2}\right)=\{0,1\} \ , \quad \cJ\left(j,\frac{1}{2},\frac{1}{2} \right)=\left\{j-1,j,j+1\right\} \ , \\
&\cM\left(j,\frac{1}{2},\frac{1}{2} \right)=\{\chi_{j+1}-\chi_{j+\frac{1}{2}},\chi_{j}-\chi_{j+\frac{1}{2}},\chi_{j+1}-\chi_{j-\frac{1}{2}},\chi_{j-1}-\chi_{j-\frac{1}{2}}\} \equiv  \{m_1,m_2,m_3,m_4\} \ ,\\
&\cM\left(\frac{1}{2},\frac{1}{2},j\right)=\{\chi_{j+1}-\chi_{1},\chi_{j}-\chi_{1},\chi_{j}-\chi_0,\chi_{j-1}-\chi_{1}\} \ .
\end{align*}
The degeneracies are $d_{j-1}=d_{j+1}=1$ and $d_{j}=2$, and the dimension of the centralizer is $\dim(\cC_{j,\frac{1}{2},\frac{1}{2}})=6$. The central elements $\alpha_i$ take the values $\alpha_1=\chi_j$ and $\alpha_2=\alpha_3=\chi_{\frac{1}{2}}$ in the quotient $\overline{AW}(j,\frac{1}{2},\frac{1}{2})$. As in the previous sections, we perform the transformation $X=(q-\qi)^2 \tilde{X}+q+\qi$ on the generators $X=A,B,D,D',K$ (see Remark \ref{rem:qlim}).
The defining relations of $\overline{AW}(j,\frac{1}{2},\frac{1}{2})$ can be written as follows 
\begin{align}
&\tB\tA\tB=\left( \tchi_j-[2]_q[1/2]_q^2\right)\tB+\tK\tB \ , \label{eq:AW1j} \\ 
&\tA\tB\tA=a_1\{\tA,\tB\}-a_2\tB+a_3\tA+\tK(\tA-a_1)+a_4 \ , \label{eq:AW2j} \\
&(\tA-\tchi_{j-\frac{1}{2}})(\tA-\tchi_{j+\frac{1}{2}})=0
\ , \quad \tB(\tB-[2]_q)=0 \ , \quad (\tK-\tchi_{j-1})(\tK-\tchi_{j})(\tK-\tchi_{j+1})=0 \ , \label{eq:ABKj} \\ 
&(\tD-\tchi_{j-\frac{1}{2}})(\tD-\tchi_{j+\frac{1}{2}})=0 \ , \quad (\tD'-\tchi_{j-\frac{1}{2}})(\tD'-\chi_{j+\frac{1}{2}})=0 \ , \label{eq:relDj}\\
&(\tK-\tB-\tchi_{j+1}+\tchi_{1})(\tK-\tB-\tchi_{j}+\tchi_{1})(\tK-\tB-\tchi_{j})(\tK-\tB-\tchi_{j-1}+\tchi_{1})=0 \ , \\
&\prod_{i=1}^{4}(\tK-\tA-m_i)=0 \ , \quad \prod_{i=1}^{4}(\tK-\tD-m_i)=0 \ , \quad \prod_{i=1}^{4}(\tK-\tD'-m_i)=0 \ , \label{eq:KADj}\\
&\Omega=\chi_{j}^2+2\chi_{\frac{1}{2}}^2+K^2+\chi_{j}\chi_{\frac{1}{2}}^2K-\chi_0^2 \ , \label{eq:Omj}
\end{align}
where
\begin{align}
\tD&=\frac{(q^2+q^{-2})}{[2]_q}(\tK+\tchi_j)+2\tchi_{\frac{1}{2}}-\frac{(q-\qi)}{[2]_q}[\tA,\tB]_q-\tA-\tB \ , \label{eq:Dj}\\
\tD'&=\frac{(q^2+q^{-2})}{[2]_q}(\tK+\tchi_j)+2\tchi_{\frac{1}{2}}-\frac{(q-\qi)}{[2]_q}[\tB,\tA]_q-\tA-\tB \ , \label{eq:D'j} \\
\Omega &= q(A+D) \chi_{\frac{1}{2}}(\chi_j+K) +q^{-1} B (\chi_{\frac{1}{2}}^2+\chi_j K) -q^2 A^2 - q^{-2} B^2 -q^2 D^2 - q ABD \ , \label{eq:Omj2}
\end{align}
and where we have used the following constants
\begin{align*}
&a_1=\frac{[2]_q[j-\frac{1}{2}]_q[j+\frac{3}{2}]_q}{q^2+q^{-2}} \ , \quad a_2=2\frac{\tchi_{j-\frac{1}{2}}\tchi_{j+\frac{1}{2}}}{q^2+q^{-2}} \ , \\
&a_3=\frac{[2]_q}{q^2+q^{-2}}(2\tchi_{\frac{1}{2}}-[2]_q[j+\textstyle\frac{1}{2}]_q^2)+\tchi_j \ , \quad a_4=a_1([j+\textstyle\frac{1}{2}]_q^2+\tchi_{\frac{1}{2}}) \ , \\
&m_1=\tchi_{j+1}-\tchi_{j+\frac{1}{2}} \ , \quad m_2=\tchi_{j}-\tchi_{j+\frac{1}{2}} \ , \quad m_3=\tchi_{j+1}-\tchi_{j-\frac{1}{2}} \ , \quad m_4=\tchi_{j-1}-\tchi_{j-\frac{1}{2}} \ .
\end{align*}
\begin{prop}\label{pr:AWj}
	The quotient $\overline{AW}(j,\frac{1}{2},\frac{1}{2})$ can be presented with the following relations
	\begin{align}
	&\tA^2=(\tchi_{j-\frac{1}{2}}+\tchi_{j+\frac{1}{2}})\tA-\tchi_{j-\frac{1}{2}}\tchi_{j+\frac{1}{2}} \ , \quad \tB^2=[2]_q\tB \ , \label{eq:ABj} \\
	&\tB\tA\tB=[2]_q\{\tA,\tB\}-[2]_q^2\tA-([j+\textstyle\frac{3}{2}]_q^2+[j-\textstyle\frac{1}{2}]_q^2-1)(\tB-[2]_q) \ . \label{eq:BABj}
	\end{align}
\end{prop}
\proof
We first show that the relations \eqref{eq:AW1j}--\eqref{eq:KADj} imply the relations of the proposition. The two equations in \eqref{eq:ABj} follow directly from \eqref{eq:ABKj}. We also deduce from the first relation of \eqref{eq:ABKj} that 
\begin{equation}
\left(\tA-a_1\right)\left(\tA-\tchi_{j-\frac{1}{2}}-\tchi_{j+\frac{1}{2}}+a_1\right)=\frac{[2j-1]_q[2j+3]_q}{q^2+q^{-2}} \ .
\end{equation}
Since the R.H.S. of the previous relation does not vanish for $j>\frac{1}{2}$, it can be used in \eqref{eq:AW2j} to find
\begin{equation}
\tK=\{\tA,\tB\} - (\tchi_{j-\frac{1}{2}}+\tchi_{j+\frac{1}{2}})\tB-[2]_q \tA+(1+[2]_q)([2j+\textstyle\frac{3}{2}]_q[\frac{1}{2}]_q+\tchi_{j-\frac{1}{2}}) \ . \label{eq:Kj} \\
\end{equation} 
Using \eqref{eq:Kj} in \eqref{eq:Dj} and \eqref{eq:D'j}, and then substituting in \eqref{eq:relDj}, one obtains expressions for $\tA\tB\tA\tB$ and $\tB\tA\tB\tA$ in terms of the elements $1,\tA,\tB,\tA\tB,\tB\tA,\tA\tB\tA$ and $\tB\tA\tB$. By using \eqref{eq:Kj} and the expressions for $\tA\tB\tA\tB$ and $\tB\tA\tB\tA$ in the third relation of \eqref{eq:ABKj}, one gets \eqref{eq:BABj}. 

Finally, we want to show that \eqref{eq:ABj} and \eqref{eq:BABj} imply the defining relations of $\overline{AW}(j,\frac{1}{2},\frac{1}{2})$. To that end, we suppose that the relations of the proposition are true and we define the element $\tK$ as in \eqref{eq:Kj}. It is then straightforward to verify that $\tK$ is central and that \eqref{eq:AW1j}--\eqref{eq:Omj} hold.
\endproof 
\begin{thm}
	Conjecture 4.1 is verified for $j_1=j$ and $j_2=j_3=\frac{3}{2}$, where $j=1,\frac{3}{2},2,...$
\end{thm}
\proof
We already know that the map $\overline{\phi}$ is surjective. From the previous proposition, we conclude that $\{1,\tA,\tB,\tA\tB,\tB\tA,\tA\tB\tA\}$ is a generating set for $\overline{AW}(j,\frac{1}{2},\frac{1}{2})$. Therefore, $\dim\left(\overline{AW}(j,\frac{1}{2},\frac{1}{2}) \right)\leq\dim \left( \cC_{j,\frac{1}{2},\frac{1}{2}} \right)=6$, which shows the injectivity of $\overline{\phi}$.
\endproof
\subsection{Connection with the one-boundary Temperley--Lieb algebra}\label{ssec:1bTL} 
On the basis of the findings for the limit $q\to1$ \cite{CPV}, one might expect that the centralizer in the case of one spin-$j$ and two spin-$\frac{1}{2}$ will be isomorphic to the one-boundary Temperley--Lieb algebra. We can indeed confirm that this algebra is recovered as a quotient of $AW(3)$.   
\begin{definition}\cite{MS,MW,NRG}
	The one-boundary Temperley--Lieb algebra $1bTL_2(q,\omega)$ is generated by $\sigma_0$ and $\sigma_1$ with the following defining relations
	\begin{equation}
	\sigma_0^2=\frac{[\omega]_q}{[\omega-1]_q}\sigma_0 \ , \quad \sigma_1^2=(q+q^{-1})\sigma_1 \ , \quad \sigma_1\sigma_0\sigma_1=\sigma_1 \ . \label{eq:1TL} 
	\end{equation}
\end{definition}
\begin{thm}
	The quotiented Askey--Wilson algebra $\overline{AW}(j,\frac{1}{2},\frac{1}{2})$, for $j=1,\frac{3}{2},2...$, is isomorphic to the one-boundary Temperley--Lieb algebra $1bTL_2(q,2j+1)$. This isomorphism is given explicitly by
	\begin{align}
	\overline{AW}\left(j,\frac{1}{2},\frac{1}{2}\right) &\to 1bTL_2(q,2j+1) \nonumber \\
	\tA&\mapsto\tchi_{j+\frac{1}{2}}-[2j]_q\sigma_0 \ , \label{eq:iso1TL1} \\
	\tB &\mapsto[2]_q-\sigma_1 \ . \label{eq:iso1TL2} 
	\end{align}
\end{thm}
\proof
It is easy to see that the map $\varphi$ is bijective. The homomorphism can be directly verified from the relations of the proposition \ref{pr:AWj}.
\endproof

\section{Conclusion and perspectives}

Summing up, we have offered a conjecture according to which a quotient of the Askey--Wilson algebra is isomorphic to the centralizer of the image of the diagonal embedding of $U_q(\gsl_2)$ in the tensor product of any three irreducible representations. It has been proved in several cases, and we thus obtained the Temperley--Lieb, Birman--Murakami--Wenzl and one-boundary Temperley--Lieb algebras as quotients of the Askey--Wilson algebra. In the limit $q \to 1$, the results of the paper \cite{CPV} are recovered. We have provided further evidence in support of the conjecture by studying the finite irreducible representations of the quotient of the Askey--Wilson algebra, more particularly in the case of three identical spins.   

Proving the conjecture in the case of three arbitrary spins $j_1,j_2,j_3$ would be an obvious continuation of the work presented here. If true, this conjecture would provide a presentation of the centralizer of $U_q(\gsl_2)$ in terms of generators and relations for any three irreducible representations. 

We could first consider, more simply, the case of three identical spins $j_1=j_2=j_3=s$. As for the Temperley--Lieb ($s=\frac{1}{2}$) and the Birman--Murakami--Wenzl ($s=1$) algebras, we expect that the centralizer for any spin $s$ will be linked to a quotient of the braid group algebra.

In \cite{CP}, a diagrammatic description of the centralizers of $U_q(gl_n)$ has been proposed. It is based on the notion of fused
Hecke algebras. Developing a connection between this diagrammatic approach and the
Askey--Wilson algebra could prove fruitful.

Throughout the present paper, we assume $q$ to be not a root of unity. This choice allows to decompose the tensor product of irreducible representations of $U_q(\gsl_2)$ into a direct sum of irreducible representations (see Subsection \ref{ssec:irreps}). As a consequence, the matrices $C_i,C_{ij},C_{123}$ are diagonalizable and their minimal polynomials are those discussed in Subsection \ref{ssec:centr}. It could be interesting to study the centralizer when $q$ is a root of unity and to examine how the quotient of $AW(3)$ is affected.

Another generalization of the results presented here would be to consider the $n$-fold tensor product of irreducible representations of $U_q(\gsl_2)$ and to connect the centralizer to a higher rank Askey--Wilson algebra $AW(n)$. The approach using the $R$-matrix proposed in \cite{CGVZ} should be helpful for this purpose. In fact, a simpler starting point could be to generalize either the conjecture given in \cite{CPV} by studying the connection between the centralizers of $\gsl_2$ and a higher rank Racah algebra, or the one given in \cite{CFV} by examining how centralizers of $\mathfrak{osp}(1|2)$ relate to the higher rank Bannai-Ito algebra $BI(n)$ (see \cite{CVZ}).

Yet another direction to generalize the results of this paper would be to study the centralizer of the diagonal embedding of $\mathfrak{g}$ or
$U_q(\mathfrak{g})$ with $\mathfrak{g}$ a higher rank Lie algebra.
A first step in this direction was made recently in \cite{CPV2} where the
centralizer $Z_2(sl_3)$ of the diagonal embedding of $sl_3$ in the
twofold tensor product of $sl_3$ has been identified.
A proposition similar to Proposition \ref{pr:surj} has also been proved in
that case. A quotient of the algebra $Z_2(sl_3)$ that describes the
centralizer for any
representations of $sl_3$ has still to be investigated. We hope to report on some of this issues in the future.

\bigskip 

{\bf Acknowledgments:}
The authors are grateful to Lo\"{\i}c Poulain D'Andecy for numerous enlightening discussions. They also thank Paul Terwilliger for useful exchanges. N. Cramp\'e is partially supported by Agence Nationale de la Recherche
Projet AHA ANR-18-CE40-0001. The work of L. Vinet is funded in part by a discovery grant of the Natural Sciences and Engineering Research Council (NSERC) of Canada. M. Zaimi holds graduate scholarships from the NSERC and the Fonds de recherche du Qu\'ebec - Nature et technologies (FRQNT).

\end{document}